\documentclass[aop,MSNbibl,seceqn,nameyear,dvips]{arximspdf}

\doi{10.1214/13-AOP837} %kopijuoti is PTS
\volume{42}
\issue{2}
\pubyear{2014}
\firstpage{623}
\lastpage{648}

\makeatletter

\newproclaim{Notation}{Notation}

\newcommand{\rrvert}{\vert}
\newcommand{\llvert}{\vert}

\newcommand{\xrightarrow}[1]{\stackrel{#1}{\rightarrow}}

\newtheorem{lemma}{Lemma}[section]
\newtheorem{proposition}{Proposition}[section]
\newtheorem{theorem}{Theorem}[section]

\newproclaim{remark}{Remark}[section]

\newcommand{\x}{\mathbf{x}}
\newcommand{\e}{\mathbb{E}}

\makeatother

\begin{document}
\begin{frontmatter}

\title{Necessary and sufficient conditions for the~asymptotic distributions
of coherence of~ultra-high dimensional random matrices}
\runtitle{Laws of coherence of random matrices}

\begin{aug}
\author[A]{\fnms{Qi-Man} \snm{Shao}\corref{}\thanksref{t1}\ead[label=e1]{qmshao@cuhk.edu.hk}}
\and
\author[B]{\fnms{Wen-Xin} \snm{Zhou}\ead[label=e2]{zhouwenxin1986@gmail.com}}
\runauthor{Q.-M. Shao and W.-X. Zhou}
\affiliation{Chinese University of Hong Kong and
Hong Kong University\break of  Science and Technology}
\address[A]{Department of Statistics\\
Chinese University of Hong Kong\\
Shatin, NT\\
Hong Kong\\
\printead{e1}}
\address[B]{Department of Mathematics\\
Hong Kong University of Science\\
\quad and Technology\\
Clear Water Bay\\
Kowloon, Hong Kong\\
\printead{e2}}
\end{aug}

\thankstext{t1}{Supported in part by Hong Kong
RGC-GRF---UST603710, CUHK2130344.}

% HISTORY:
\received{\smonth{1} \syear{2012}}
\revised{\smonth{12} \syear{2012}}

% ABSTRACT
%
\begin{abstract}
Let $\x_1,\ldots, \x_n$ be a random sample from a $p$-dimensional
population distribution, where $p=p_n \to\infty$ and $\log p =
o(n^\beta)$ for some $0 <\beta\leq1$, and let $L_n$ be the
\textit{coherence} of the sample correlation matrix. In this paper it
is proved that $\sqrt{n/\log p} L_n \to2 $ in probability if and
only if $Ee^{t_0 |x_{11}|^\alpha} < \infty$ for some $t_0>0$, where
$\alpha$ satisfies $ \beta= \alpha/(4-\alpha)$. Asymptotic
distributions of $L_n$ are also proved under the same sufficient
condition. Similar results remain valid for $m$-\textit{coherence} when
the variables of the population are $m$ dependent. The proofs are
based on self-normalized moderate deviations, the Stein--Chen method
and a newly developed randomized concentration inequality.
\end{abstract}

% KEYWORDS
\begin{keyword}[class=AMS]
\kwd{60F05}
\kwd{62E20}
\end{keyword}
\begin{keyword}
\kwd{Coherence}
\kwd{sample correlation matrix}
\kwd{extreme distribution}
\kwd{law of large numbers}
\end{keyword}

\end{frontmatter}

%s1 #&#
\section{Introduction}\label{sec1}

This paper is motivated by the recent results of Cai and Jiang (\citeyear{CJ11a},
\citeyear{CJ11b}) on asymptotic behaviors of the largest magnitude of off-diagonal
entries of the sample correlation matrix. Consider a $p$-variable
population represented by a random vector $\mathbf{x}=(x_1,\ldots,
x_p)^T$ with the covariance matrix $\Sigma$, and let $X_n=(x_{ij})$ be
an $n\times p$ random matrix where the $n$ rows consist a random sample
of size $n$ from the population. The Pearson correlation coefficient
$\rho_{ij}$ between the $i$th and $j$th columns of $X_n$ is given by
%
%e1.1 #&#
%
\begin{equation}\label{rho}
\rho_{ij} = \frac{\sum_{k=1}^n(x_{ki}-\bar{x}_i)(x_{kj}-\bar
{x}_j)}{\sqrt{\sum_{k=1}^n
(x_{ki}-\bar{x}_i)^2} \cdot\sqrt{\sum_{k=1}^n(x_{kj}-\bar
{x}_j)^2}},\qquad 1\leq i, j \leq p,
\end{equation}
where $\bar{x}_i=(1/n)\sum_{k=1}^n x_{ki}$. Then the sample
correlation matrix
$\Gamma_n$ is defined by $\Gamma_n \equiv(\rho_{ij})$.

The main object of interest in this paper is the largest magnitude of
off-diagonal entries of the sample correlation matrix, that is,
%
%e1.2 #&#
%
\begin{equation}\label{coherence}
L_n = \max_{1\leq i<j\leq p}|\rho_{ij}|.
\end{equation}
As in \citet{CJ11a}, $L_n$ is called the \textit
{coherence} of the random matrix~$X_n$.

In the case where $p$ and $n$ are of the same order, that is, $n/p
\rightarrow\lambda\in(0, \infty)$, asymptotic properties of coherence
$L_n$ have been extensively studied recently. \citet{J04} was the
first to establish the strong laws and limiting distributions of $L_n$.
The moment assumption in \citet{J04}
has been substantially improved by \citet{LR06}, \citet{Wa07},
\citet{LLS08}, \citet{LLR09} and \citet{LQR10}. \citet{LLS08} proved that similar results
hold for $p=O(n^{\alpha})$ where $\alpha$ is a constant. We refer to
\citet{CJ11a} and references therein for recent developments on
this topic. In particular, \citet{CJ11a} considered the
ultra-high dimensional case where $p$ can be as large as $e^{n^\beta}$
for some $\beta\in(0,1)$. Specifically, assuming all the entries of
$X_n$, $\{x_{ij}, i\geq1, j\geq1 \}$ are i.i.d. real-valued random
variables with mean $\mu$ and variance $0<\sigma^2< \infty$, they
proved the following results.\vspace*{9pt}

\textit{Suppose $\e e^{t_0 |x_{11}|^{\alpha}}<\infty$ for some $t_0>0$
and $\alpha>0$. Assume that $p =p_n\rightarrow\infty$ and $\log
p=o(n^{\beta})$ as $n\rightarrow\infty$}, \textit{where $\beta=\frac{\alpha
}{4+\alpha}$. Then}
%
%e1.3 #&#
%
\begin{equation}\label{weaklaw}
\sqrt{n/(\log p)} L_n \rightarrow2\qquad \mbox{in probability.}
\end{equation}
\textit{If $0<\alpha\leq2$}, \textit{then}
%
%e1.4 #&#
%
\begin{equation}\label{limitlaw}
nL_n^2 - 4\log p+ \log_2 p \xrightarrow{d.}
Y,
\end{equation}
\textit{where $d$. denotes convergence in distribution}, \textit{$\log_2 p \equiv\log
\log p$ and the random variable $Y$ has an extreme distribution of type
I with distribution function}
%
%e1.5 #&#
%
\begin{eqnarray}\label{extre}
F_Y(y)&=& e^{-(1/\sqrt{8\pi})e^{-y/2}},\qquad
y\in\mathbb{R}.\\[-14pt]\nonumber
\end{eqnarray}

The main purpose of this paper is to find necessary and sufficient
conditions for (\ref{weaklaw}) and (\ref{limitlaw}). Our result shows that
the optimal choice of $\beta$ is that
$
\beta= \alpha/(4-\alpha)$, $ 0 < \alpha\leq2$
for (\ref{weaklaw}), and the same $\beta$ for (\ref{limitlaw}) when $0<
\alpha\leq1$. It is also shown that, when $1<\alpha\leq4/3$ and $\e
(x_{11}-\mu)^3 \neq0$, (\ref{limitlaw}) does not hold, but a recentered
$L_n$ will do.

The rest of the paper is organized as follows. The main results,
Theorems \ref{thm1}, \ref{thm2} and \ref{thm3} will be stated in
Section~\ref{sec2}.
A closely related problem of testing for $m$-dependence of the
population is
considered and an application to compressed sensing is revisited in
this section.
The proofs of Theorems \ref{thm1} and \ref{thm2} are given in Sections~\ref{sec3} and \ref{sec4}, respectively,
by using the Stein--Chen method, moderate deviations for both
standardized and self-normalized
sums of independent random variables.
The proof of Theorem \ref{thm3} is postponed to Section~\ref{sec5}.

%s2 #&#
\section{Main results}\label{sec2}

In this section, we consider the law of large
numbers and asymptotic distributions of the \textit{coherence} $L_n$.
In Section~\ref{sec2.1}, we provide necessary and sufficient conditions for the
two aforementioned limiting properties and the optimal choice of $\beta
$ in terms of $\alpha$. In Section~\ref{sec2.2}, we consider the $m$-\textit
{coherence}, $L_{n,m}$, of a random matrix with $m$-dependent structure
in each row.

\begin{Notation*}
Throughout this paper, $a_n \asymp b_n$ will denote
that there exist two positive constants $c_1$, $c_2$ such that $c_1
\leq a_n/b_n \leq c_2$, for all $n\geq1$; $a_n \sim b_n$ will denote
$\lim_{n\rightarrow\infty} a_n/b_n =1$.
\end{Notation*}

%s2.1 #&#
\subsection{The i.i.d. case}\label{sec2.1} \label{case1} In this subsection, we
assume that the entries $x_{ij}$ of $X_n$ are i.i.d. with mean $\mu$
and variance $\sigma^2>0$. Let
%
%e2.1 #&#
%
\begin{equation}\label{alpha}
\beta=\beta_{\alpha}=\alpha/(4-\alpha),\qquad 0< \alpha\leq2.
\end{equation}

We first state the law of large numbers for $L_n$.

%th2.1 #&#
\begin{theorem} \label{thm1}
\textup{(i)} Suppose $\e\exp\{ t_0|x_{11}|^{\alpha} \}< \infty
$ for some $0<\alpha\leq2$ and $t_0>0$. Assume $p=p_n\rightarrow
\infty$ and $\log p=o(n^{\beta_{\alpha}}) $ as $n \rightarrow\infty$. Then
%
%e2.2 #&#
%
\begin{equation}\label{wl}
\sqrt{n/( \log p )} L_n \rightarrow2
\end{equation}
in probability as $n \rightarrow\infty$.

\textup{(ii)} Let $ 0 < \beta\leq1$. If (\ref{wl}) holds for any
$p\to\infty$ satisfying
$\log p =o(n^{\beta})$, then $\e\exp\{ t_0|x_{11}|^{\alpha}\} <
\infty$ for some $t_0>0$, where $\alpha=\alpha_{\beta}=4\beta
/(1+\beta)$; that is, $\alpha$ and $\beta$ satisfy (\ref{alpha}).
\end{theorem}

%re2.1 #&#
\begin{remark}
Clearly, when $\alpha=2$, $\beta$ equals to 1, so the range for
dimension $p$ reduces to $\log p = o(n)$. On the other hand, as proved
by \citet{CJ11b}, if $x_{11} \sim\mathcal{N}(0,1)$ and $(\log
p)/n \rightarrow\gamma\in(0, \infty)$, then
\[
L_n \rightarrow\sqrt{1-e^{-4\gamma}}>0 \qquad\mbox{in probability
as } n \rightarrow\infty.
\]
Hence, result (\ref{wl}) no longer holds for $\log p \asymp n$. We
believe that the limit of $L_n$ will also depend on the distribution of
$x_{11}$ in this case, which still remains an open question.
\end{remark}

The next theorem gives the asymptotic distribution of $L_n$ after
proper normalization. Let $\kappa= \e(x_{11}-\mu)^3/\sigma^3$ and
%
%e2.3 #&#
%
\begin{equation}
\label{wn} W_n = \cases{ nL_n^2-4\log p+
\log_2 p, &\quad $0<\alpha\leq1$, \vspace*{2pt}
\cr
nL_n^2-4
\log p - \bigl(8\kappa^2/3\bigr)n^{-1/2}(\log
p)^{3/2}\cr
\qquad{}+\log_2 p, &\quad $1< \alpha\leq4/3$.}
\end{equation}

%th2.2 #&#
\begin{theorem} \label{thm2}
Suppose\vspace*{1pt} $\e\exp\{ t_0|x_{11}|^{\alpha} \} < \infty$ for some $0<
\alpha\leq4/3$ and \mbox{$t_0>0$}. Assume $p=p(n)\rightarrow\infty$, $\log
p=o(n^{\beta_{\alpha}}) $ as $n \rightarrow\infty$.
Then
%
%e2.4 #&#
%
\begin{equation}\label{re2}
W_n \stackrel{d.} {\to} Y,
\end{equation}
where $Y$ has the distribution function given in (\ref{extre}).
\end{theorem}

Clearly, when $\alpha=4/3$, $\beta_\alpha= 1/2$, (\ref{re2}) converges
weakly to the distribution function (\ref{extre}) provided that
$\log p =o(n^{1/2})$. However, (\ref{re2}) is not valid when $\log p
\asymp n^{1/2}$ as shown in \citet{CJ11b}; that is, if $x_{11}
\sim\mathcal{N}(0,1)$ and $(\log p) /n^{1/2} \rightarrow\gamma\in
[0,\infty)$, the limiting\vspace*{1pt} distribution of (\ref{limitlaw}) is shifted to
the left by $8\gamma^2$, that is, $\exp\{-(1/\sqrt{8\pi})
e^{-(y+8\gamma
^2)/2}\}$, $y \in\mathbb{R}$. For $ 4/3 < \alpha\leq2$, derivation
of the limiting distribution of $L_n$ needs more delicate arguments.

Theorems \ref{thm1} and \ref{thm2} together fully exhibit the
dependence between ranges of dimension $p$ and the optimal moment
conditions for asymptotic properties (\ref{weaklaw}) and (\ref{limitlaw})
of the coherence $L_n$.

%re2.2 #&#
\begin{remark}\label{rmk2}
It is known that the convergence rate to type I extreme distribution is
typically slow. When $p \asymp n$, \citet{LLS08} proved that
the rate of convergence can be improved to $O ((\log
n)^{5/2}n^{-1/2}  )$ if an ``intermediate'' approximation is used,
that is,
%
%e2.5 #&#
%
\begin{eqnarray}\label{inter}
&&\sup_{y\in\mathbb{R}} \biggl| P\bigl(nL_n^2
\leq y\bigr)
- \exp \biggl\{ -\frac{p(p-1)}{2}P\bigl(\chi^2_1
\geq y\bigr) \biggr\} \biggr| \nonumber\\[-8pt]\\[-8pt]
&&\qquad= O \biggl( \frac{(\log n)^{5/2}}{n^{1/2}} \biggr),
\nonumber
\end{eqnarray}
where $\chi^2_1$ has a chi-square distribution with one degree of
freedom. In the ultra-high dimensional case, Theorem \ref{thm2} implies
%
%e2.6 #&#
%
\begin{eqnarray}\label{inter1}\quad
&&\sup_{y\in\mathbb{R}} \biggl| P(W_n \leq y)
- \exp \biggl\{ -\frac{p(p-1)}{2}P\bigl(\chi^2_1
\geq4\log p-\log\log p+y\bigr) \biggr\} \biggr| \nonumber\\[-8pt]\\[-8pt]
&&\qquad\rightarrow0.\nonumber
\end{eqnarray}
It is possible to prove that the rate of convergence of (\ref{inter1}) is
of order $O(n^{-1/2})$. To test the independence of the $p$-variate
population, it may be better to choose the critical value based on the
``intermediate'' approximation. That is, reject the null hypothesis if
$L_n^2 \geq z_{\alpha}/n$, where $z_{\alpha}$ satisfies $P(\chi_1^2
\geq z_{\alpha})=-2\log(1-\alpha)/\allowbreak\{p(p-1)\}$.
\end{remark}

%re2.3 #&#
\begin{remark}\label{rmk3}
Both Theorems \ref{thm1} and \ref{thm2} are still valid if $L_n$ is
replaced by
%
%e2.7 #&#
%
\begin{equation}\label{cs-cohere}
\tilde{L}_n =\max_{1\leq i< j \leq p} |\tilde{
\rho}_{ij}|,
\end{equation}
where
%
%e2.8 #&#
%
\begin{equation}\label{cs-rho1}
\tilde{\rho}_{ij} = \frac{\sum_{k=1}^n(x_{ki}-\mu)(x_{kj}-\mu
)}{\sqrt{ \sum_{k=1}^n
(x_{ki}-\mu)^2 \sum_{k=1}^n(x_{kj}-\mu)^2 }}.
\end{equation}
The quantity $\tilde{L}_n$ arises from compress sensing literature.
See, for example, \citet{Do06}.
\end{remark}

%s2.2 #&#
\subsection{$m$-dependent case}\label{sec2.2} \label{case2}
As discussed in \citet{CJ11a}, a variant of coherence $L_n$ can
be used to construct a test for bandedness of the covariance matrix in
the Gaussian case. In this paper, we drop the normality assumption and
consider a more general problem of testing whether the population is
$m$-dependent, where $m$ can depend on $n$. More specifically, let
$X_n= (x_{ij})_{n\times p}$, where the $n$ rows are i.i.d. random
vectors drawn from a $p$-variate population represented by $\mathbf
{x}=(x_1,\ldots, x_p)^T$ with the covariance matrix $\Sigma$. Assume
all $p$ components of $\mathbf{x}$ are identically distributed with
mean $\mu$ and variance $\sigma^2>0$. Then, we wish to test the hypothesis
%
%e2.9 #&#
%
\begin{equation}\label{band}
H_0\dvtx  x_i \mbox{ and } x_j \mbox{ are
independent for all } |i-j|\geq m.
\end{equation}

Analogous to the definition of $L_n$, we introduce the \textit
{$m$-coherence} of the matrix $X_n$ as follows:
%
%e2.10 #&#
%
\begin{equation}\label{banded}
L_{n,m}= \max_{|i-j|\geq m} |\rho_{ij}|.
\end{equation}
In addition, let $(r_{ij} )_{p \times p}$ be the correlation matrix of
$\mathbf{x}$.
For any given $0<\delta<1$, set
%
%e2.11 #&#
%
\begin{equation}\label{Ga}\quad
\Gamma_{p,\delta} = \bigl\{1\leq i\leq p\dvtx  |r_{ij}|>1-\delta\mbox{ for
some } 1\leq j\leq p \mbox{ with } j\neq i \bigr\}.
\end{equation}
The following theorem establishes the limiting distribution of
$L_{n,m}$ under the null hypothesis.

%th2.3 #&#
\begin{theorem} \label{thm3}
Let $\kappa=\e(x_{11}-\mu)^3/\sigma^3$ and define
\[
W_{n,m} = \cases{ nL_{n,m}^2-4\log p+
\log_2 p, &\quad $0<\alpha\leq1$, \vspace*{2pt}
\cr
nL_{n,m}^2-4
\log p - \bigl(8\kappa^2/3\bigr)n^{-1/2}(\log
p)^{3/2}+\log_2 p, &\quad $1< \alpha\leq4/3$.}
\]
Suppose $\e\exp\{ t_0|x_{11}|^{\alpha} \} < \infty$ for some $0<
\alpha\leq4/3$ and $t_0>0$.
Moreover, assume that, as $n\rightarrow\infty$:
\begin{longlist}[(ii)]
\item[(i)] $p=p_n \rightarrow\infty$, $\log p=o(n^{\beta
_{\alpha}}) $, where $\beta_{\alpha}$ is given in (\ref{alpha});
\item[(ii)] there exists some $\delta\in(0,1) $ such that
$\llvert  \Gamma_{p,\delta} \rrvert  = o(p)$
and $m=o(p^{\varepsilon_\delta})$, where $\varepsilon_\delta= ( 2 \delta
- \delta^2)/(4-2\delta+\delta^2)$.
\end{longlist}
Then, under $H_0$, $W_{n,m}$ converges weakly to the extreme
distribution (\ref{extre}).
\end{theorem}

Theorem \ref{thm3} was proved in \citet{CJ11a} when $\mathbf
{x}$ is multivariate normal, $\log p=o(n^{1/3})$, $m=o(p^t)$ for any
$t>0$ and $| \Gamma_{p,\delta} |=o(p)$ for some $\delta\in(0,1)$.
It was also pointed out therein that the assumption $| \Gamma
_{p,\delta} |=o(p) $ is essential in the sense that there exists a
covariance matrix $\Sigma$ such that the conclusion of Theorem~\ref
{thm3} for Gaussian entries no longer holds when $p \sim n
e^{n^{1/4}}$, $m=n$ and $\llvert  \Gamma_{p,\delta} \rrvert =p$ for
any $\delta>0$.
In Theorem \ref{thm3} here, the assumption on $m$ is weakened, and
condition (i) provides the optimal choice of $\beta$ in terms of
$\alpha$, and more importantly, Gaussian entries are not required.

%re2.4 #&#
\begin{remark}
Similar to Remark \ref{rmk2}, an ``intermediate'' approximation can
also be applied here based on
%
%e2.12 #&#
%
\begin{eqnarray}\label{inter2}
&&\sup_{y\in\mathbb{R}} \bigl| P(W_{n,m} \leq y)
\nonumber\\[-8pt]\\[-8pt]
&&\qquad{} - \exp \bigl\{ -\bigl(p^2/2\bigr)P\bigl(\chi^2_1
\geq4\log p-\log\log p+y\bigr) \bigr\} \bigr| \rightarrow0
\nonumber
\end{eqnarray}
as $n \rightarrow\infty$.
\end{remark}

%re2.5 #&#
\begin{remark}
In compressed sensing, the quantity $\tilde{L}_n$, defined in (\ref
{cs-cohere}), is useful because it is closely related to the so-called
\textit{mutual incoherence property} (MIP), which requires the
pairwise correlations among column vectors of $X=X_{n\times p}$ to be
small. More specifically, under certain assumptions on $X$, the condition
%
%e2.13 #&#
%
\begin{equation}
\label{MIP}
(2k-1)\tilde{L}_n <1
\end{equation}
guarantees the exact recovery of $\beta\in\mathbb{R}^p$ from linear
measurements $y=X \beta$, when $\beta$ has at most $k$ nonzero
entries. This condition is also sharp in the sense that there exists
matrices $X_0$ such that recovering some $k$-sparse signals $\beta$
based on $y=X_0\beta$ when $(2k-1)\tilde{L}_n=1$ is impossible. See,
\citet{Do01}, \citet{F04} and \citet{CJ10}.

It was shown in \citet{CJ11a} that the limiting properties of $\tilde
{L}_n$ can be directly applied to compute the probability that random
measurement matrices satisfy the MIP conditions (\ref{MIP}). In
particular, Theorem \ref{thm1} with $L_n$ replaced with $\tilde{L}_n$
provides necessary and sufficient conditions for
$ \tilde{L}_n \sim2\sqrt{(\log p)/n} $. This suggests that the sparsity
$k$ should satisfy
$
k< \sqrt{n/(\log p)}/4
$
approximately in order for the MIP condition (\ref{MIP}) to hold.
\end{remark}

%s3 #&#
\section{\texorpdfstring{Proof of Theorem \protect\ref{thm1}}{Proof of Theorem 2.1}}\label{sec3}

We start with collecting some technical lemmas that will be used to
prove our main results. Without loss of generality, assume $\{x_{ij};
1\leq i\leq n, 1\leq j\leq p \}$ are i.i.d. random variables with mean
zero and variance one. Both letters $C$ and $c$ denote constants that
do not depend on $n$ or $p$, but may depend on the distribution of
$x_{11}$ and vary from line to line.

%s3.1 #&#
\subsection{Technical lemmas}\label{sec3.1} \label{techtool}
As in many previous works on the extreme distribution approximation,
the following lemma is a special case of Theorem 1 of
\citet{Ar89}, based on the Stein--Chen method.

%le3.1 #&#
\begin{lemma} \label{pa}
Let $\{\eta_{\alpha}, \alpha\in I\}$ be random variables on an index
set I. For each $\alpha\in I$, let $B_\alpha$ be a subset of I with
$\alpha\in B_\alpha$. For any given $t \in\mathbb{R}$, set
$\lambda=\sum_{\alpha\in I} P(\eta_\alpha>t)$. Then
%
%e3.1 #&#
%
\begin{equation}\label{pa1}
\Bigl| P \Bigl( \max_{\alpha\in I} \eta_\alpha\leq t \Bigr)
-e^{-\lambda
} \Bigr| \leq\min\bigl(1, \lambda^{-1}\bigr)
(b_{1}+b_{2}+b_{3}),
\end{equation}
where
\begin{eqnarray*}
b_{1}&=&\sum_{\alpha\in I}\sum
_{\beta\in B_\alpha}P(\eta_\alpha >t)P(\eta_{\beta}>t),\qquad
b_{2}=\sum_{\alpha\in I}\mathop{\sum
_{\beta\in B_\alpha}}_{\beta\neq
\alpha} P(\eta_\alpha>t, \eta_{\beta}>t),
\\
b_{3} &=& \sum_{\alpha\in I} \e\bigl|P\bigl(
\eta_\alpha>t|\sigma(\eta_\beta, \beta\notin B_\alpha)
\bigr) -P(\eta_\alpha>t)\bigr|
\end{eqnarray*}
and $\sigma(\eta_\beta, \beta\notin B_\alpha)$ is the $\sigma
$-algebra generated
by $\{\eta_\beta, \beta\notin B_\alpha\}$. In particular, if $\eta
_\alpha$ is
independent of $\{\eta_\beta, \beta\notin B_\alpha\}$, for each
$\alpha\in I$, then $b_{3}$ vanishes.
\end{lemma}

For a sequence of random variables $X_1, X_2, \ldots\,$, we use $S_n$ and
$V_n^2$ to denote the partial sum and the partial quadratic sum,
respectively, that is,
\[
S_n = \sum_{i=1}^n
X_i,\qquad V_n^2=\sum
_{i=1}^n X_i^2.
\]
The following lemma is due to \citet{LY61} on the moderate deviation
under i.i.d. assumption.
%
%le3.2 #&#
\begin{lemma} \label{ld}
Suppose $X_1, X_2, \ldots$ are i.i.d. random variables with $\e X_1 =0$
and $\e X_1^2=1$:
\begin{longlist}[(iii)]
\item[(i)] If $\e e^{t_0|X_1|^\alpha}< \infty$ for some
$0<\alpha\leq1$ and $t_0>0$, then
%
%e3.2 #&#
%
\begin{equation}\label{ldp}
\lim_{n \rightarrow\infty} \frac{1}{x_n^2}\log P ( S_n/\sqrt
{n} \geq x_n ) = -1/2
\end{equation}
for any $x_n \rightarrow\infty$, $x_n=o ( n^{{\alpha
}/({2(2-\alpha)})}  )$.

\item[(ii)] If $\e e^{t_0|X_1|^\alpha}< \infty$ for some
$0<\alpha\leq1/2$ and $t_0>0$, then
%
%e3.3 #&#
%
\begin{equation}\label{mdp1}
\frac{P(S_n/\sqrt{n} \geq x)}{1-\Phi(x)} \rightarrow1
\end{equation}
holds uniformly for $0\leq x \leq o ( n^{{\alpha}/({2(2-\alpha
)})}  )$.
\item[(iii)] Assume $\e e^{t_0 X_1}< \infty$ for some $t_0>0$.
If $x\geq0$, $x=o(n^{1/4})$, then
%
%e3.4 #&#
%
\begin{equation}\label{mdp}
\frac{P(S_n/\sqrt{n} \geq x)}{1-\Phi(x)} = \exp \biggl\{\frac{x^3\e
X_1^3}{6n^{1/2}} \biggr\} \biggl[1+ O \biggl(
\frac{1+x}{n^{1/2}} \biggr) \biggr].\vadjust{\goodbreak}
\end{equation}
\end{longlist}
\end{lemma}

We also need the following self-normalized moderate deviations:

%le3.3 #&#
\begin{lemma}[{[\citet{Shao97}]}] \label{smd}
Assume that $X_1, X_2, \ldots$ are i.i.d. random variables with $\e X_1
=0$ and $0< \sigma^2= \e X_1^2 < \infty$. Then, for any sequence of
real numbers $x_n$ satisfying $x_n \rightarrow\infty$ and
$x_n=o(\sqrt{n})$,
%
%e3.5 #&#
%
\begin{equation}\label{sld}
\log P(S_n/ V_n \geq x_n)
\thicksim-x_n^2/2.
\end{equation}
\end{lemma}

%s3.2 #&#
\subsection{\texorpdfstring{Proof of Theorem \protect\ref{thm1}}{Proof of Theorem 2.1}}\label{sec3.2} \label{pf1}

\mbox{}

\begin{pf*}{Proof of \textup{(i)}}
The main idea of the proof is to show that $L_n$ can be reduced to
$L_{n,0}=\max_{1 \leq i < j \leq p} |\rho_{ij,0}|$, where
%
%e3.6 #&#
%
\begin{equation}\label{cs-rho}
\rho_{ij,0} = \frac{1}{n\sigma^2}\sum_{k=1}^n(x_{ki}-
\mu ) (x_{kj}-\mu),\qquad 1\leq i, j\leq p.
\end{equation}
Let
%
%e3.7 #&#
%
\begin{eqnarray}\label{not1}
S_{n,i} &=& \sum_{k=1}^n
x_{ki},\qquad V_{n,i}^2 =\sum
_{k=1}^n x_{ki}^2,\nonumber\\[-8pt]\\[-8pt]
\Delta_{n,i}&=&\frac
{S_{n,i}}{\sqrt{n}V_{n,i}},\qquad
1\leq i\leq p, n\geq1.\nonumber
\end{eqnarray}
Decompose the sample correlation coefficient as
%
%e3.8 #&#
%
\begin{equation}\label{decom-rho}
\rho_{ij} = \rho_{ij,1}-\rho_{ij,2},\qquad 1\leq i, j
\leq p
\end{equation}
and accordingly, define
\[
L_{n,k} =\max_{1\leq i<j\leq p} |\rho_{ij,k}|,\qquad k=1,
2,
\]
where
%
%e3.9 #&#
%
\begin{eqnarray}\label{rho12}
\rho_{ij,1} &=& \frac{\sum_{k=1}^n x_{ki}x_{kj}/(V_{n,i}V_{n,j}) }{
\{(1-\Delta_{n,i}^2)(1-\Delta_{n,j}^2) \}^{1/2}},\nonumber\\[-8pt]\\[-8pt]
\rho_{ij,2} &=&
\frac{\Delta_{n,i}\Delta_{n,j} }{
\{(1-\Delta_{n,i}^2)(1-\Delta_{n,j}^2) \}^{1/2}}.\nonumber
\end{eqnarray}
Intuitively, Lemma \ref{smd} suggests that $\Delta_{n,i}$ can be
negligible and Lemma \ref{ld} indicates that $V_{n,i}^2/n$ is close to
$1$. Let
%
%e3.10 #&#
%
\begin{equation}\label{en12}
\varepsilon_{n1} = c_1 (\log p)^{1/2}/n^{\beta/2}
\quad\mbox{and}\quad \varepsilon_{n2} = c_2 (\log
p)^{1/2}/n^{1/2},
\end{equation}
where $c_1$ and $c_2$ are positive constants only depending on the
distribution of $x_{11}$ and will be specified later in different
cases. Since
$
\e\exp\{ t_0|x_{11}^2-1|^{\alpha/2} \}< \infty
$, it follows from (\ref{ldp}) and (\ref{sld}) that
%
%e3.11 #&#
%
\begin{equation}\label{ld1}
P \bigl( \bigl|V_{n,1}^2-n\bigr| / n^{1/2} >
\varepsilon_{n1} n^{\beta/2} \bigr) \leq2\exp\bigl\{ -c
\varepsilon_{n1}^2 n^{\beta} \bigr\}
\end{equation}
and
%
%e3.12 #&#
%
\begin{equation}\label{smld1}
P\bigl( |\Delta_{n,1}| > \varepsilon_{n2} \bigr) \leq2\exp\bigl\{ - c
\varepsilon _{n2}^2 n \bigr\}
\end{equation}
for all sufficiently large $n$. Now define the subset
%
%e3.13 #&#
%
\begin{equation}\label{event1}
\mathcal{E}_n = \Bigl\{ \max_{1\leq i\leq p}
\bigl|V_{n,i}^2/n-1\bigr| \leq \varepsilon_{n1}n^{(\beta-1)/2},
\max_{1\leq i\leq p} |\Delta_{n,i}| \leq\varepsilon_{n2}
\Bigr\}.
\end{equation}
Then, for properly chosen $c_1$ and $c_2$ in (\ref{en12}), we have
%
%e3.14 #&#
%
\begin{equation}\label{pr-event1}
P\bigl(\mathcal{E}_n^c\bigr) \leq2p \bigl( \exp\bigl\{-c
\varepsilon_{n1}^2 n^{\beta}\bigr\} +\exp\bigl\{-c
\varepsilon _{n2}^2 n \bigr\} \bigr) = o\bigl(p^{-4}
\bigr).
\end{equation}
Recall $L_{n,0}$ defined through (\ref{cs-rho}). Clearly, on
$\mathcal{E}_n$
\[
\frac{L_{n,0}}{1+\varepsilon_{n1} n^{(\beta-1)/2}} \leq L_{n,1} \leq \frac{L_{n,0}}{(1-\varepsilon_{n2}^2)(1-\varepsilon_{n1} n^{(\beta-1)/2})}
\]
and
\[
L_{n,2} \leq\varepsilon_{n2}^2/\bigl(1-
\varepsilon_{n2}^2\bigr).
\]
Noting that $\varepsilon_{n1}n^{(\beta-1)/2}= c_1(\log p)^{1/2}/n^{1/2}
= o(1)$ and
$\sqrt{ n/\log p} \varepsilon_{n2}^2 =\break  c_2^2 (\log p)^{1/2}/n^{1/2} = o(1)$,
we have on $\mathcal{E}_n$
%
%e3.15 #&#
%
\begin{equation}\label{Ln0-1}
L_{n,1}/L_{n,0} \to1,\qquad \sqrt{n /\log p} | L_n
- L_{n,1}| \to0,
\end{equation}
which together with (\ref{pr-event1}) shows that conclusion (\ref
{wl}) will be
a direct consequence of the next proposition. The proof is postponed to
the end of this section.
\end{pf*}

%pr3.1 #&#
\begin{proposition} \label{prop1}
Under the conditions of \textup{(i)} in Theorem \ref{thm1}, we have
$
\sqrt{n/(\log p)} L_{n,0} \rightarrow2$ in probability as
$n\rightarrow
\infty$.
\end{proposition}

\begin{pf*}{Proof of \textup{(ii)}}
We shall prove the necessity of moment conditions under a weaker
assumption than (\ref{wl}). Assume that there exists a constant $C_0
\geq4$, such that
%
%e3.16 #&#
%
\begin{equation}\label{nec1}
P \Bigl( \sqrt{n/(\log p)} \max_{1\leq i<j \leq p}|\rho_{ij}|
\geq C_0 \Bigr) \rightarrow0.
\end{equation}
Note that
$
\max_{1\leq i< j\leq p} |\rho_{ij}| \geq\max_{1\leq i \leq p/2} |
\rho_{i,[p/2]+i}| $, then (\ref{nec1}) implies
%
%e3.17 #&#
%
\begin{equation}\label{nec2}
P \Bigl( \max_{1\leq i\leq p/2} |\rho_{i, [p/2]+i}| > C_0
\sqrt{(\log p)/n} \Bigr) \rightarrow0.
\end{equation}
Observe that $\{ \rho_{i, [p/2]+i}$, $1\leq i\leq[p/2]\}$ are i.i.d.
random variables and
that $\sum_{k=1}^n (x_{ki} -\bar{x}_i)^2 \leq\sum_{k=1}^n
x_{ki}^2$, (\ref{nec2}) thus yields
%
%e3.18 #&#
%
\begin{equation}\label{nec3}
p \cdot P \biggl( \frac{|\sum_{k=1}^n x_{k1}x_{k2}- n \bar{x}_1 \bar
{x}_2|}{(\sum_{k=1}^n
x_{k1}^2)^{1/2} (\sum_{k=1}^n x_{k2}^2)^{1/2}} > C_0\sqrt{(\log p)/n} \biggr)
\rightarrow0.
\end{equation}
For $n\geq16$, define the subset
\[
\mathcal{D}_n = \biggl\{ \frac{\sum_{k=2}^n x_{ki}^2}{n} \leq2,
{
|\sum_{k=2}^n x_{ki}| \over\sqrt{n}} \leq n^{1/4}, i=1, 2;
\frac{ |\sum_{k=2}^n x_{k1}x_{k2}|}{\sqrt{n}}
\leq1 \biggr\}.
\]
By the central limit theorem and the strong law of large numbers,
$P(\mathcal{D}_n) \rightarrow2\Phi(1)-1$,
so that $P(\mathcal{D}_n) \geq1/2$ for sufficiently large $n$.
Furthermore, since $\log p =o(n)$, we have on $\mathcal{D}_n$,
\begin{eqnarray*}
&&\biggl\{ \frac{|\sum_{k=1}^n x_{k1}x_{k2}|}{(\sum_{k=1}^n x_{k1}^2)^{1/2}
(\sum_{k=1}^n x_{k2}^2)^{1/2}} > C_0\sqrt{
\frac{\log p}{n}} \biggr\}
\\
&&\qquad\supseteq \biggl\{ \frac{ |x_{11}x_{12}| - 2 \sqrt{n} - |x_{11}| - |x_{12}|
}{(x_{11}^2+2n)^{1/2}(x_{12}^2+2n)^{1/2}} > C_0\sqrt{
\frac{\log p}{n}} \biggr\}
\\
&&\qquad\supseteq \bigl\{ \bigl(|x_{11}|-c \sqrt{\log p}\bigr) \bigl(|x_{12}| -
c \sqrt{\log p} \bigr) > 3C_0\sqrt{n \log p} \bigr\}
\end{eqnarray*}
for some $c>0$, which along with the independence of $\mathcal{D}_n$
and $\{x_{11}, x_{12}\}$ yields
%
%e3.19 #&#
%
\begin{eqnarray}\label{nec3-0}\quad
&&
P \biggl( \frac{|\sum_{k=1}^n x_{k1}x_{k2} - n \bar{x}_1 \bar{x}_2
|}{(\sum_{k=1}^n
x_{k1}^2)^{1/2} (\sum_{k=1}^n x_{k2}^2)^{1/2}} > C_0 \sqrt{(\log p)/n}
\biggr)
\nonumber\\
&&\qquad\geq P(\mathcal{D}_n) \cdot P \bigl( \bigl(|x_{11}|-c\sqrt{
\log p}\bigr) \bigl(|x_{12}| - c\sqrt{\log p} \bigr) > 3C_0\sqrt{n \log
p} \bigr)
\\
&&\qquad\geq (1/2) \cdot \bigl\{ P \bigl(|x_{11}| > 2C_0^{1/2}(n
\log p)^{1/4} \bigr) \bigr\}^2.
\nonumber
\end{eqnarray}
If follows from (\ref{nec3}) and (\ref{nec3-0}) that
%
%e3.20 #&#
%
\begin{equation}\label{nec4}
p^{1/2} P \bigl(|x_{11}| > C_0 (n \log
p)^{1/4}\bigr) = o(1)
\end{equation}
for any\vspace*{1pt} $p$ satisfying $\log p = o(n^{\beta})$. By a contradiction
argument, it is easy to see that
(\ref{nec4}) implies that
$\e\exp\{ t_0|x_{11}|^{4\beta/(1+\beta)} \} < \infty$, for some
$t_0>0$. This proves part (ii).
\end{pf*}

We end this section with the proof of Proposition \ref{prop1}.

%s3.3 #&#
\subsection{\texorpdfstring{Proof of Proposition \protect\ref{prop1}}{Proof of Proposition 3.1}}\label{sec3.3}

It suffices to show, for any $0<\varepsilon< 1/8$, as
$n\rightarrow\infty$,
%
%e3.21 #&#
%
\begin{equation}\label{wl1}
P \bigl( \sqrt{n/(\log p)}L_{n,0} \leq2-\varepsilon \bigr) \rightarrow0
\end{equation}
and
%
%e3.22 #&#
%
\begin{equation}\label{wl2}
P \bigl( \sqrt{n/(\log p)}L_{n,0} > 2 + \varepsilon \bigr) \rightarrow0.
\end{equation}
We apply Lemma \ref{pa} to prove (\ref{wl1}) by using (\ref{pa1}) to deal
with the maximum.
% then appealing to large deviations for the probability of each
%individual event.
%New techniques are needed to estimate joint probabilities.
The proof of (\ref{wl2}) is similar, and so the details are omitted here.

Put $y_n=(2- \varepsilon)\sqrt{(\log p) /n}$, $n\geq1$. Define
\[
I=\bigl\{(i,j); 1\leq i<j\leq p\bigr\},\qquad
A_{ij} = \bigl\{ |
\rho_{ij,0}| > y_n \bigr\},\qquad 1\leq i<j \leq p,
\]
and
\[
B_{i,j}= \bigl\{ (k,l)\in I \setminus\bigl\{(i, j)\bigr\}
\mbox{; either } k \in \{i, j\} \mbox{ or } l \in\{i, j\} \bigr\}.
\]
Since $\{ x_{ij}; (i,j)\in I \}$ are identically distributed, by Lemma
\ref{pa},
%
%e3.23 #&#
%
\begin{equation}\label{0t}
\Bigl| P \Bigl( \max_{1\leq i< j\leq p} |\rho_{ij,0}| \leq (2-
\varepsilon)\sqrt{(\log p)/n} \Bigr) - e^{-\lambda_n} \Bigr| \leq b_{n,1} +
b_{n,2},
\end{equation}
where
%
%e3.24 #&#
%
\begin{eqnarray}\label{3t}
\lambda_n&=&\frac{p(p-1)}{2}P(A_{12}),\qquad
b_{n,1}\leq p^3P^2(A_{12}),\nonumber\\[-8pt]\\[-8pt]
b_{n,2}&\leq& p^3P(A_{12}A_{13}).\nonumber
\end{eqnarray}
Because $0< \alpha/2 \leq1$ and $\e\exp\{t_0|x_{11}x_{12}|^{\alpha
/2}\}<\infty$,
it follows from (\ref{ldp}) that, for all sufficiently large $n$,
%
%e3.25 #&#
%
\begin{eqnarray}\label{P12}
P(A_{12}) &=& P \biggl( \frac{|\sum_{k=1}^n x_{k1}x_{k2} | }{n^{1/2}} > \sqrt
{n}y_n \biggr)
\nonumber\\[-8pt]\\[-8pt]
&\leq& 2\exp\bigl\{-(1-\varepsilon)ny_n^2/2\bigr\} =
2p^{-(1-\varepsilon
)(2-\varepsilon)^2/2},\nonumber
\end{eqnarray}
which, in turn implies
%
%e3.26 #&#
%
\begin{equation}\label{2t}
\lambda_n \rightarrow\infty \quad\mbox{and}\quad b_{n,1}=o(1)
\qquad\mbox{as } n \rightarrow\infty.
\end{equation}
As for $b_{n,2}$, we have
%
%e3.27 #&#
%
\begin{eqnarray}\label{b2-1}
P(A_{12}A_{13}) & = & P \biggl( \frac{|\sum_{k=1}^n x_{k1}x_{k2}|}{n} >
y_n, \frac{|\sum_{k=1}^n
x_{k1}x_{k3}|}{n} > y_n \biggr)
\nonumber\\
&\leq& P \biggl( \frac{| \sum_{k=1}^n x_{k1}(x_{k2}+x_{k3}) |}{n} > 2y_n \biggr)
\\
&&{} + P \biggl( \frac{| \sum_{k=1}^n x_{k1}(x_{k2}- x_{k3}) |}{n} > 2y_n \biggr).
\nonumber
\end{eqnarray}
Since $\e[x_{k1}(x_{k2}+x_{k3})]=0$ and $\e
[x_{k1}(x_{k2}+x_{k3})]^2=2$, applying (\ref{ldp}) again, we get
\[
P \biggl( \frac{| \sum_{k=1}^n x_{k1}(x_{k2}+x_{k3}) |}{n} > 2y_n \biggr) \leq2\exp\bigl\{-(1-
\varepsilon)ny_n^2\bigr\} = 2p^{-(1-\varepsilon)(2-\varepsilon)^2}.
\]
Similarly, the same result holds for $P(|\sum_{k=1}^n x_{k1}(x_{k2}-
x_{k3})| >
2 y_n n)$.
Therefore,
%
%e3.28 #&#
%
\begin{equation}\label{b2-3}
b_{n,2}\leq p^3P(A_{12}A_{13})=O
\bigl( p^{3-(1-\varepsilon)(2-\varepsilon
)^2} \bigr) = o(1).
\end{equation}
This completes the proof of (\ref{wl1}) by (\ref{0t}), (\ref{3t}),
(\ref{2t})
and (\ref{b2-3}).

%s4 #&#
\section{\texorpdfstring{Proof of Theorem \protect\ref{thm2}}{Proof of Theorem 2.2}}\label{sec4}

The main idea is to use Lemma \ref{pa} again.
The proof of part (i) is standard while that of part (ii)
requires a more delicate estimate of
$\lambda_n$ given in (\ref{3t}). In particular, we need a randomized
concentration inequality
in Lemma \ref{rci}.

We formulate the proof into two cases.\vspace*{9pt}

\textit{Case} 1. $0 < \alpha\leq1$.\vspace*{9pt}

For arbitrary fixed $y\in\mathbb{R}$, let
%
%e4.1 #&#
%
\begin{equation}\label{zn-0}
y_n=\sqrt{(y+4\log p-\log_2 p )/n},\qquad
\log_2 p \equiv\log\log p
\end{equation}
for large $n$ so that $y+4\log p-\log_2 p >0$. We need to prove that
%
%e4.2 #&#
%
\begin{equation}\label{cp1}
P \Bigl( \max_{1\leq i<j\leq p} |\rho_{ij}| \leq y_n
\Bigr) \rightarrow\exp \bigl(-(1/\sqrt{8\pi})e^{-z/2} \bigr).
\end{equation}
Similar to (\ref{0t}), we have
%
%e4.3 #&#
%
\begin{equation}\label{ineq1}
\Bigl| P \Bigl( \max_{1\leq i<j\leq p} |\rho_{ij}| \leq
y_n \Bigr) - e^{-\lambda_n} \Bigr| \leq b_{n,1}+b_{n,2},
\end{equation}
where $\lambda_n$, $b_{n,1}$, $b_{n,2}$ and $A_{ij}$ are defined as in
(\ref{3t}) with $\rho_{ij,0}$ replaced by $\rho_{ij}$. It suffices
to show
%
%e4.4 #&#
%
\begin{equation}\label{PA12}
P(A_{12}) \sim2 \bigl( 1-\Phi(\sqrt{n}y_n) \bigr)+o
\bigl(p^{-2}\bigr) \sim \frac
{e^{-y/2}}{\sqrt{2\pi}}p^{-2}
\end{equation}
and
%
%e4.5 #&#
%
\begin{equation}\label{PA123}
P(A_{12}A_{13}) = o\bigl(p^{-3}\bigr).
\end{equation}
Analogously to (\ref{event1}), let
%
%e4.6 #&#
%
\begin{equation}\label{event2}
\mathcal{E}_{n\cdot3} = \Bigl\{ \max_{i=1, 2, 3}
\bigl|V_{n,i}^2/n-1\bigr| \leq\varepsilon_{n1}n^{(\beta-1)/2},
\max_{i=1, 2, 3} |\Delta_{n,i}| \leq
\varepsilon_{n2} \Bigr\},
\end{equation}
where $V_{n,i}$ and $\Delta_{n,i}$ are given in (\ref{not1}). In view of
(\ref{pr-event1}),
we can choose $c_1$ and $c_2$ in (\ref{en12}) properly such that
%
%e4.7 #&#
%
\begin{equation}\label{00-1}
P\bigl(\mathcal{E}_{n\cdot3}^c\bigr)=o\bigl(p^{-3}
\bigr).
\end{equation}
On $\mathcal{E}_{n\cdot3}$, we have
%
%e4.8 #&#
%
\begin{equation}\label{r21u}
|\rho_{1i}| \leq\frac{|\rho_{1i,0}|}{(1-\varepsilon
_{n2}^2)(1-\varepsilon
_{n1}n^{(\beta-1)/2})} +
\frac{\varepsilon_{n2}^2}{1-\varepsilon_{n2}^2},\qquad
i=2,3,
\end{equation}
and [recall $y_n \sim2n^{-1/2}(\log p)^{1/2}$]
%
%e4.9 #&#
%
\begin{equation}\label{r21a}
|\rho_{12}|= \bigl\{ 1+o \bigl( \sqrt{(\log p)/n} \bigr) \bigr\}
\cdot |\rho_{12,0}| + O \bigl((\log p)/n \bigr).
\end{equation}

We are now ready to prove (\ref{PA12}) and (\ref{PA123}).

\begin{pf*}{Proof of (\ref{PA12})}
By (\ref{r21a}), it follows that, on $\mathcal{E}_{n\cdot3}$,
\[
\bigl\{ |\rho_{12}|>y_n \bigr\} = \bigl\{ |\rho_{12,0}| >
\hat{y}_n \bigr\}\qquad \mbox{with } \hat{y}_n =
y_n \bigl( 1+o\bigl(n^{-1/2}(\log p)^{1/2}\bigr)
\bigr).
\]
Recalling the definition of $\rho_{12,0}$ in (\ref{cs-rho}) and
\[
\e x_{k1}x_{k2}=0,\qquad \e(x_{k1}x_{k2})^2=1,\qquad
\e e^{t_0|x_{11}x_{12}|^{\alpha/2}}< \infty \qquad\mbox{with } 0< \alpha /2 \leq1,
\]
it follows directly from (\ref{mdp1}) that, as $n\rightarrow\infty$,
%
%e4.10 #&#
%
\begin{equation}\label{md-r12}
\frac{ P( \rho_{12,0} >\hat{y}_n ) }{1-\Phi(\sqrt{n}\hat{y}_n)} \rightarrow1.
\end{equation}
Noticing that $\log p =o(n^{1/3})$, it is easy to check that
\[
\frac{1-\Phi(\sqrt{n}y_n)}{1-\Phi(\sqrt{n}\hat{y}_n)} \rightarrow1,
\]
which, together with (\ref{md-r12}) yields (\ref{PA12}).
\end{pf*}

\begin{pf*}{Proof of (\ref{PA123})}
By (\ref{r21u}), following the same argument as in (\ref{b2-1}) and
(\ref{b2-3}), we have for any $0< \varepsilon< 1/8$,
\begin{eqnarray*}
&&P(A_{12}A_{13})
\\
&&\qquad\leq P \bigl( |\rho_{12,0}| \geq\bigl\{1-o(1)\bigr\}y_n,
|\rho_{13,0}| \geq \bigl\{1-o(1)\bigr\}y_n \bigr) +P\bigl(
\mathcal{E}_{n\cdot3}^c\bigr)
\\
&&\qquad\leq C \exp\bigl\{-(1-\varepsilon)ny_n^2\bigr\} +o
\bigl(p^{-3}\bigr)
\\
&&\qquad\leq C (\log p) p^{-4(1-\varepsilon)} + o\bigl(p^{-3}\bigr) = o
\bigl(p^{-3}\bigr).
\end{eqnarray*}
This gives (\ref{PA123}).\vspace*{9pt}

%%\subsection{Proof of Theorem \ref{thm2},

\textit{Case} 2. $ 1 < \alpha\leq4/3$.\vspace*{9pt}

Similar to $y_n$ in (\ref{zn-0}), for $y\in\mathbb{R}$ we now define
%
%e4.11 #&#
%
\begin{equation}\label{newzn}
y_n=\sqrt{(y+4\log p+c_{n,p}-\log_2 p) /n},
\end{equation}
where $c_{n,p}=(8\kappa^2/3)n^{-1/2}(\log p)^{3/2}$. Following the
same argument as in the proof of case 1, (\ref{PA123}) remains valid. It
thus remains to show that
%
%e4.12 #&#
%
\begin{equation}\label{newp12}
P(A_{12}) \sim2\mathcal{L}_{n,y}+o\bigl(p^{-2}
\bigr),
\end{equation}
where
\[
\mathcal{L}_{n,y} = \bigl(1-\Phi(\sqrt{n}y_n) \bigr)\exp
\bigl(\kappa^2ny_n^3/6\bigr).
\]
Let $\x_i = ( x_{i1},\ldots, x_{ni})^T$, $i=1,\ldots, p$ be the $p$
columns of $X_n$, and $\| \cdot\|$ denotes the Euclidean norm in
$\mathbb{R}^n$.
Rewrite $\rho_{12}$ as
%
%e4.13 #&#
%
\begin{eqnarray}\label{not2}
\rho_{12} = \hat{\rho}_{12}/ \bigl\{\bigl(1-
\Delta_{n,1}^2\bigr) \bigl(1-\Delta _{n,2}^2
\bigr)\bigr\} ^{1/2} \nonumber\\[-8pt]\\[-8pt]
\eqntext{\mbox{with } \displaystyle \hat{\rho}_{12}
\equiv\frac{\mathbf{x}_1^T \mathbf{x}_2-
n^{-1}S_{n,1}S_{n,2}}{
\|\mathbf{x}_1\|\|\mathbf{x}_2\|}.}
\end{eqnarray}
Define the subset
%
%e4.14 #&#
%
\begin{equation}\label{event2}
\mathcal{E}_{n\cdot2}=\bigl\{\max\bigl(|\Delta_{n,1}|, |
\Delta_{n,2}|\bigr) \leq \varepsilon_{n2}\bigr\},
\end{equation}
where $\varepsilon_{n2} = c_2(\log p)^{1/2}/n^{1/2}$ is given in (\ref
{en12}) with $c_2>0$ chosen appropriately such that $P(\mathcal
{E}_{n\cdot2}^c)=o(p^{-4})$. Hence, with probability at least $1-o(p^{-4})$,
%
%e4.15 #&#
%
\begin{equation}\label{r12app}
|\rho_{12}|/|\hat{\rho}_{12}| = 1+o\bigl(n^{-1/2}
\bigr).
\end{equation}

For $\hat{\rho}_{12}$, using the elementary inequalities
\[
2ab \leq a^2+b^2 \quad\mbox{and}\quad
(1+s)^{1/2} \geq1+s/2-s^2/2\qquad \mbox{for any } s>-1
\]
to give lower and upper bounds as follows:
%
%e4.16 #&#
%
\begin{equation}\label{r12lb}
\{ \hat{\rho}_{12} > y_n \} \supseteq \bigl\{
\mathbf{x}_1^T \mathbf{x}_2 - y_n
\bigl(\| \mathbf{x}_1\|^2+\| \mathbf{x}_2
\|^2\bigr)/2 > n^{-1}S_{n,1}S_{n,2}
\bigr\}
\end{equation}
and
%
%e4.17 #&#
%
\begin{eqnarray}\label{r12ub}
&&\{ \hat{\rho}_{12} > y_n \}
\nonumber\\
&&\qquad \subseteq \bigl\{ \mathbf{x}_1^T
\mathbf{x}_2 - y_n \bigl(\| \mathbf{x}_1
\|^2+\| \mathbf{x}_2\|^2\bigr)/2
\\
&&\hspace*{15.6pt}\qquad > n^{-1}S_{n,1}S_{n,2} - ny_n^2
\bigl[ \bigl(\|\mathbf{x}_1\|^2/n-1\bigr)^2 +
\bigl(\| \mathbf{x}_2\|^2/n-1\bigr)^2 \bigr]
\bigr\}.
\nonumber
\end{eqnarray}
Therefore, in order to prove (\ref{newp12}), we need to show the
following two claims:
%
%e4.18 #&#
%
\begin{equation}\label{asym1}
P \bigl( \mathbf{x}_1^T \mathbf{x}_2 -
y_n \bigl(\| \mathbf{x}_1\|^2+\|
\mathbf{x}_2\|^2\bigr)/2 > 0 \bigr) \sim
\mathcal{L}_{n,y}+o\bigl(p^{-2}\bigr)
\end{equation}
and
%
%e4.19 #&#
%
\begin{equation}\label{asym2}\quad
P \bigl( \Delta_n < \mathbf{x}_1^T
\mathbf{x}_2 - y_n \bigl(\| \mathbf {x}_1\|
^2+\| \mathbf{x}_2\|^2\bigr)/2 \leq0 \bigr) =
o(1) \bigl\{\mathcal{L}_{n,y}+p^{-2} \bigr\},
\end{equation}
where $\Delta_n=\Delta(S_{n,1}, S_{n,2}, V_{n,1}^2, V_{n,2}^2)$ is
given by
%
%e4.20 #&#
%
\begin{equation}\label{Den1}
\Delta_n= n^{-1}S_{n,1}S_{n,2} -
ny_n^2 \bigl[ \bigl(\|\mathbf{x}_1
\|^2/n-1\bigr)^2 + \bigl(\|\mathbf{x}_2
\|^2/n-1\bigr)^2 \bigr].
\end{equation}
\upqed\end{pf*}

\begin{pf*}{Proof of (\ref{asym1})} Given two random vectors $\mathbf
{x}_1, \mathbf{x}_2 \in\mathbb{R}^n$, truncate one of which as follows:
%
%e4.21 #&#
%
\begin{equation}\label{trun}
x_{k2}^{\tau} = x_{k2}I_{\{|x_{k2}|\leq\tau\}},\qquad k=1,\ldots,
n, \qquad\mbox{with } \tau= \tau_n =
t_0^{-1/\alpha}n^{\beta/\alpha}\hspace*{-15pt}
\end{equation}
and write
%
%e4.22 #&#
%
\begin{equation}\label{xi}
\xi_k= \xi_{n,k} = y_nx_{k1}
x_{k2}^{\tau}-y_n^2
\bigl(x_{k1}^2+x_{k2}^{\tau2}\bigr)/2,\qquad
k=1,\ldots, n.
\end{equation}
By the union bound and Markov inequality,
%
%e4.23 #&#
%
\begin{equation}\label{tail}
P \Bigl( \max_{1\leq k\leq n } |x_{k2}| > \tau \Bigr) \leq\e
\bigl[e^{t_0|x_{11}|^{\alpha}} \bigr] \cdot ne^{-n^{\beta}}
\end{equation}
and it is easy to see that $\mathbf{x}_1^T \mathbf{x}_2 - y_n (\|
\mathbf{x}_1\|^2+\| \mathbf{x}_2\|^2)/2 =y_n^{-1}\sum_{k=1}^n\xi_k
$ on $\{
\max_{k}|x_{k2}|\leq\tau\}$. We thus aim to estimate the probability
$P(\sum_{k=1}^n\xi_k >0)$. Since $\alpha>1$ and
$
y_n \tau^{2-\alpha} = O ((\log p)^{1/2}/n^{\beta/2}  ) = o(1)
$, it follows that
\begin{eqnarray*}
\xi_k &\leq& y_n\tau^{2-\alpha}|x_{k1}||x_{k2}|^{\alpha-1}
\leq y_n\tau^{2-\alpha} \bigl(|x_{k1}|^{\alpha}
+ |x_{k2}|^{\alpha
} \bigr) \\
&=&o(1) \bigl(|x_{k1}|^{\alpha}+|x_{k2}|^{\alpha}
\bigr),
\end{eqnarray*}
which, in turn, implies $\sup_{1\leq k\leq n, n\geq1} \e e^{\xi_k} <
\infty$. Moreover, it is easy to verify
that
\begin{eqnarray*}
\e\xi_k &=& -y_n^2 + y_n^2
\e x_{11}^2I_{\{|x_{11}|>\tau\}}/2 = -y_n^2
\bigl\{1+O\bigl(y_n^2\bigr) \bigr\},
\\
\operatorname{Var}(\xi_k) &=& y_n^2 \bigl\{1+O
\bigl(y_n^2\bigr) \bigr\} \quad\mbox{and}\quad
\frac{\e(\xi_k-\e\xi_k)^3}{\operatorname{Var}^{3/2}(\xi_k)} = \bigl(\e x_{11}^3\bigr)^2 +
O(y_n).
\end{eqnarray*}
Let $\mu_n=\sum_{k=1}^n\e\xi_k$ and $\sigma_n^2=\sum_{k=1}^n
\operatorname{Var}(\xi_k)$, then
$-\mu_n/\sigma_n=\sqrt{n}y_n\{1+O(y_n^2)\}$. Moreover, noting that
$\sqrt {n}y_n=o(n^{1/4})$ and $\kappa= \e x_{11}^3$ (with $\mu=0$ and
$\sigma
^2=1$), it follows from (\ref{mdp}) and the above facts that
\begin{eqnarray*}
P \Biggl( \sum_{k=1}^n
\xi_k >0 \Biggr) &=& P \biggl( \frac
{\sum_{k=1}^n(\xi_k- \e
\xi_k) }{\sigma_n}> -\mu_n/
\sigma_n \biggr)
\\
&\sim& \bigl(1-\Phi(-\mu_n/\sigma_n) \bigr)\exp \biggl(
\frac{(-\mu
_n/\sigma_n)^3}{6n^{1/2}} \bigl(\kappa^2+O(y_n) \bigr) \biggr)
\\
&\sim& \bigl(1-\Phi(\sqrt{n}y_n) \bigr)\exp \biggl\{
\frac{\kappa
^2ny_n^3}{6} \biggr\} = \mathcal{L}_{n,y}\qquad \mbox{as } n
\rightarrow\infty.
\end{eqnarray*}
This, along with (\ref{tail}), implies (\ref{asym1}) immediately.
\end{pf*}

\begin{pf*}{Proof of (\ref{asym2})}
This requires a more delicate analysis. The main idea is to apply a
combination of the multivariate conjugate method and a randomized
concentration inequality to the truncated variables as defined in
(\ref{xi}) and (\ref{trun}). Further to the notation used in the
proof of
(\ref{asym1}), let $\{ \mathbf{y}_k= (x_{k1}, x_{k2}^{\tau});1 \leq
k\leq
n\}$ be a sequence of independent $\mathbb{R}^2$-valued random
variables and let measurable function $g\dvtx  \mathbb{R}^2
\rightarrow\mathbb{R}^3$ be given by
%
%e4.24 #&#
%
\begin{equation}\label{cj2}
\forall(u,v) \in\mathbb{R}^2\qquad g(u,v)= \bigl(uv, u^2,
v^2 \bigr).
\end{equation}
Put
\[
\mathbf{S}_n = \sum_{k=1}^n
\mathbf{y}_k= \Biggl( \sum_{k=1}^n
x_{k1}, \sum_{k=1}^n
x_{k2}^{\tau
} \Biggr)^T
\]
and
\[
\mathbf{V}_n = \sum_{k=1}^n
g(\mathbf{y}_k) = \Biggl(\sum_{k=1}^n
x_{k1}x_{k2}^{\tau}, \sum
_{k=1}^n x_{k1}^2, \sum
_{k=1}^n x_{k2}^{\tau2}
\Biggr)^T.
\]
Let\vspace*{1pt} $\lambda_n = (y_n, -y_n^2/2, -y_n^2/2)^T \in\mathbb{R}^3 $. Observe
that $\xi_k=\xi_{n,k}$ given in (\ref{xi}) can be rewritten as
$\lambda
_n^T g(\mathbf{y}_k) $ that satisfy
%
%e4.25 #&#
%
\begin{equation}\label{mgf2d}
\max_{1\leq k\leq n, n\geq1} m_{n,k} < \infty,
\end{equation}
where
\[
m_{n,k} =\e e^{\xi_k} = \e\bigl[e^{ \lambda_n^T g(\mathbf{y}_k)}\bigr].
\]
Now, let $\hat{\mathbf{y}}_1, \hat{\mathbf{y}}_2,\ldots, \hat
{\mathbf{y}}_n$ be a
sequence of independent $\mathbb{R}^2$-valued random variables such
that $\hat{\mathbf{y}}_k$ has the following distribution:
%
%e4.26 #&#
%
\begin{equation}\label{cj-dist}
\forall B \in\mathcal{B}^2\qquad
P( \hat{\mathbf{y}}_k
\in B ) = \frac
{1}{m_{n,k}} \e\bigl[e^{\lambda_n^T g(\mathbf{y}_k)}I_{\{\mathbf{y}_k \in B \}}\bigr].
\end{equation}
Accordingly, put $\hat{\mathbf{S}}_n = \sum_{k=1}^n\hat{\mathbf{y}}_k$,
$\hat{\mathbf{V}}_n = \sum_{k=1}^n g(\hat{\mathbf{y}}_k)$. The multivariate
conjugate method says that, for any $C\in\mathcal{B}^5$,
%
%e4.27 #&#
%
\begin{equation}\label{m-cj}
P\bigl\{ (\mathbf{S}_n, \mathbf{V}_n) \in C \bigr\} = \e
\bigl[e^{\lambda_n^T \hat{\mathbf{V}}_n}I_{\{(\hat{\mathbf
{S}}_n,\hat
{\mathbf{V}}_n) \in C\}}\bigr] \prod_{k=1}^n
m_{n,k}.
\end{equation}
In particular, define subsets
\begin{eqnarray*}
C_n &=& \bigl\{\mathbf{u}\in\mathbb{R}^5\dvtx
\Delta(u_1,u_2,u_4,u_5) \leq
u_3-y_n(u_4+u_5)/2 < 0 \bigr\}
\cap E_n,
\\
E_n &=& \biggl\{ \mathbf{u} \in\mathbb{R}^3 \times
\mathbb{R}^2_+\dvtx  \frac{u_1}{\sqrt{u_4}} \leq \varepsilon_{n2}n^{1/2},
\biggl|\frac{u_j}{n}-1 \biggr| \leq\varepsilon_{n1}n^{(\beta-1)/2}, j=4, 5
\biggr\},
\end{eqnarray*}
where in accordance with (\ref{Den1}),
%
%e4.28 #&#
%
\begin{equation}\label{Den2}
\Delta(v_1, v_2, v_3, v_4)=
n^{-1}v_1v_2-ny_n^2
\bigl[(v_3/n-1)^2+(v_4/n-1)^2
\bigr]
\end{equation}
and $\{ \varepsilon_{n1}, \varepsilon_{n2}; n\geq1\}$ are given as in
(\ref{en12}), such that
%
%e4.29 #&#
%
\begin{equation}\label{tail2}
P\bigl\{(\mathbf{S}_n, \mathbf{V}_n) \in
E_n^c \bigr\} =o\bigl(p^{-4}\bigr).
\end{equation}

By (\ref{m-cj}), we have
%
%e4.30 #&#
%
\begin{eqnarray}\label{m-cj-1}
P\bigl\{ (\mathbf{S}_n, \mathbf{V}_n) \in
C_n \bigr\} &=& \Biggl( \prod_{k=1}^n
m_{n,k} \Biggr) \times \e\bigl[e^{-\lambda_n^T \hat{\mathbf{V}}_n}I_{\{(\hat{\mathbf{S}}_n,
\hat
{\mathbf{V}}_n) \in C_n\} } \bigr]
\nonumber\\[-8pt]\\[-8pt]
&:= & \Biggl( \prod_{k=1}^n
m_{n,k} \Biggr) \times K_n.\nonumber
\end{eqnarray}
Let $ \hat{\xi}_k = \lambda_n^T g(\hat{\mathbf{y}}_k) $ be the
conjugate version of $\xi_k$. Then, by (\ref{cj-dist}),
\[
\e\hat{\xi}_k = \e\bigl[\xi_k e^{\xi_k}\bigr]/
\e\bigl[e^{\xi_k}\bigr],\qquad
\operatorname{Var}(\hat{\xi}_k) = \e\bigl[
\xi_k^2 e^{ \xi_k}\bigr]/ \e\bigl[e^{\xi_k}
\bigr] - (\e \hat{\xi}_k)^2.
\]
Put $\hat{\mu}_n=\sum_{k=1}^n\e\hat{\xi}_k$ and $\hat{\sigma
}_n^2 = \sum_{k=1}^n
\operatorname{Var}(\hat{\xi}_k )$.
Routine calculations show (recall $\kappa=\e x_{11}^3$)
\begin{eqnarray*}
\e\bigl[e^{ \xi_k}\bigr] &=& 1-y_n^2/2+
\kappa^2 y_n^3 /6 + O\bigl(y_n^4
\bigr),
\\
\e\bigl[\xi_k e^{ \xi_k}\bigr] &=& \kappa^2
y_n^3/2 + O\bigl(y_n^4\bigr),
\\
\e\bigl[\xi_k^2 e^{ \xi_k}\bigr] &=&
y_n^2 +\kappa^2 y_n^3
+ O\bigl(y_n^4\bigr).
\end{eqnarray*}
Consequently,
%
%e4.31 #&#
%
\begin{equation}\label{mn}
\hat{\mu}_n =\kappa^2 ny_n^3/2+
O\bigl(ny_n^4\bigr),\qquad
\hat{\sigma}_n^2=
ny_n^2+\kappa^2ny_n^3
+ O\bigl(ny_n^4\bigr)
\end{equation}
and
%
%e4.32 #&#
%
\begin{equation}\label{prod-mgf}
\prod_{k=1}^n m_{n,k} = \exp
\bigl( -ny_n^2/2 + \kappa^2ny_n^3/6+O
\bigl( ny_n^4\bigr) \bigr).
\end{equation}
As for $K_n$ in (\ref{m-cj-1}), we shall show that
%
%e4.33 #&#
%
\begin{equation}\label{kn-asym}
\sqrt{n}y_n K_{n} = o(1).
\end{equation}
Now combining (\ref{m-cj-1}), (\ref{prod-mgf}), (\ref{kn-asym}) and the
well-known result $1-\Phi(s) \sim(2\pi)^{-1/2}s^{-1}e^{-s^2/2}$ as
$s\rightarrow\infty$, it follows
\[
P\bigl\{ (\mathbf{S}_n, \mathbf{V}_n) \in
C_n \bigr\} = o(\mathcal{L}_{n,y}).
\]
This, together with (\ref{tail}), (\ref{tail2}) and the definition of
$C_n$, gives (\ref{asym2}).

\begin{pf*}{Proof of (\ref{kn-asym})}
Observe that on the event $\{ (\hat{\mathbf{S}}_n, \hat{\mathbf{V}}_n)
\in C_{n} \}$,
%
%e4.34 #&#
%
\begin{equation}\label{dehat}
\lambda_n^T \hat{\mathbf{V}}_n = \sum
_{k=1}^n\hat{\xi}_k \geq
(y_n/n)\hat {S}_{n,1}\hat{S}_{n,2} -
2n^{\beta}y_n^3\varepsilon_{n1}^2,
\end{equation}
where $ \hat{S}_{n,1} = \sum_{k=1}^n\hat{x}_{k1}$, $ \hat{S}_{n,2}
= \sum_{k=1}^n
\hat{x}_{k2}^{\tau}$. Using H\"older's inequality gives
%
%e4.35 #&#
%
\begin{eqnarray}\label{kn2}
K_{n} &\leq& \bigl( \e e^{-2\lambda_n^T \hat{\mathbf{V}}_n} I_{\{ (\hat{\mathbf{S}}_n, \hat{\mathbf{V}}_n) \in C_{n} \}}
\bigr)^{1/2}
\nonumber\\
&&{} \times \Biggl( P \Biggl( (y_n/n)\hat{S}_{n,1}
\hat{S}_{n,2} - 2n^{\beta
}y_n^3
\varepsilon_{n1}^2 \leq \sum_{k=1}^n
\hat{\xi}_k <0 \Biggr) \Biggr)^{1/2}
\\
&:= & K_{n,1}^{1/2} \times K_{n,2}^{1/2}.
\nonumber
\end{eqnarray}
We first estimate $K_{n,1}$. By (\ref{cj-dist}),
\begin{eqnarray*}
\e[\hat{x}_{k1}] &=& m_{n,k}^{-1}\e
\bigl[x_{k1}e^{ \xi_k}\bigr] = -\kappa y_n^3
/2+ O\bigl(y_n^4\bigr),
\\
\e\bigl[\hat{x}_{k1}^2\bigr] &=& m_{n,k}^{-1}
\e\bigl[x_{k1}^2e^{\xi_k}\bigr] = 1-
y_n^2/2 - \kappa^2 y_n^3/2
+ O\bigl(y_n^4\bigr)
\end{eqnarray*}
and same expansions hold for $\e[\hat{x}_{k2}^{\tau}]$ and $\e[\hat
{x}_{k2}^{\tau2}]$ as well. Thus, for all sufficiently large $n$,
$\sum_{k=1}^n\e\hat{x}_{k2}^{\tau2} \leq n $ and on $\{ (\hat
{\mathbf
{S}}_n, \hat{\mathbf{V}}_n) \in C_{n} \}$,
\[
|\hat{S}_{n,1}| \leq\sqrt{2}\varepsilon_{n2}n,\qquad
\sum_{k=1}^n\hat {x}_{k2}^{\tau2}
\leq2n.
\]
In view of (\ref{en12}) and (\ref{dehat}),
%
%e4.36 #&#
%
\begin{eqnarray}\label{dehat-ubd}\qquad
-2\lambda_n^T \hat{\mathbf{V}}_n &\leq&
-2(y_n/n)\hat{S}_{n,1} (\hat {S}_{n,2}-\e
\hat{S}_{n,2}) -2y_n\e\bigl[\hat{x}_{12}^{\tau}
\bigr]\hat {S}_{n,1} +4n^{\beta}y_n^3
\varepsilon_{n1}^2
\nonumber\\[-8pt]\\[-8pt]
&\leq& Cn^{-1/2}(\log p)Z_n + O\bigl(n^{-3/2}(\log
p)^{5/2}\bigr),\nonumber
\end{eqnarray}
where
\[
Z_n \equiv\frac{|\sum_{k=1}^n(\hat{x}_{k2}^{\tau}-\e\hat
{x}_{k2}^{\tau
})|}{4\sqrt{\sum_{k=1}^n \operatorname{Var}(\hat{x}_{k2}^{\tau})} + \sqrt{\sum_{k=1}^n(\hat
{x}_{k2}^{\tau}-\e\hat{x}_{k2}^{\tau})^2}}.
\]
Now we can use the following sub-Gaussian property of self-normalized
sums [see Lemma 6.4 in \citet{JSW03}]:

%le4.1 #&#
\begin{lemma} \label{JSW}
Let $\{ X_i, 1\leq i\leq n \}$ be a sequence of independent random
variables with $\e X_i=0$ and $\e X_i^2<\infty$. Then, for $a>0$,
\[
P \Biggl( \Biggl| \sum_{i=1}^n
X_i \Biggr| \geq a \Biggl( 4D_n + \Biggl( \sum
_{i=1}^n X_i^2
\Biggr)^{1/2} \Biggr) \Biggr) \leq8e^{-a^2/2},
\]
where $D_n^2 =\sum_{i=1}^n \e X_i^2$.
\end{lemma}

Indeed, Lemma \ref{JSW} implies $P(Z_n \geq a) \leq8e^{-a^2/2}$,
$\forall a>0$. Hence,
\[
\forall t>0\qquad
\e e^{tZ_n} \leq1+8\sqrt{2\pi}te^{t^2/2},
\]
which together with (\ref{dehat-ubd}) yields
%
%e4.37 #&#
%
\begin{equation}\label{kn21}
K_{n,1} = O(1).
\end{equation}
Next, we estimate $K_{n,2}$. The key technical tool is
the randomized concentration inequality below developed in \citet{SZ11}:

%le4.2 #&#
\begin{lemma} \label{rci}
Let $\eta_1,\ldots,\eta_n$ be independent random variables,
\[
W_n=\sum_{k=1}^n\eta_k
\]
and let $\Delta_1=\Delta_1(\eta_1,\ldots,\eta_n)$ and $\Delta_2=\Delta
_2(\eta_1,\ldots,\eta_n) $ be two measurable functions of
$\eta_1,\ldots,\eta_n$. Assume that
\[
\e\eta_k=0 \qquad\mbox{for } k=1,2, \ldots, n\quad \mbox{and}\quad \sum
_{k=1}^n\e\eta_k^2=1.
\]
For each $1\leq k \leq n$, let $\Delta_1^{(k)}$ and $\Delta_2^{(k)}$
be any
random variables such
that $\eta_k$ and $(\Delta_1^{(k)},\Delta_2^{(k)},W_n-\eta_k)$ are
independent. Then
\begin{eqnarray*}
&&P(\Delta_1 \leq W_n \leq\Delta_2)
\\
&&\qquad\leq 21 \Biggl( \sum_{k=1}^n\e|
\eta_k|^3+ \e|\Delta_2-\Delta_1|
\\
&&\hspace*{16.2pt}\qquad\quad{} + \sum_{k=1}^n\bigl\{ \e\bigl|
\eta_k\bigl(\Delta_1-\Delta_1^{(k)}
\bigr)\bigr| +\e\bigl|\eta _k\bigl(\Delta-\Delta _2^{(k)}
\bigr)\bigr| \bigr\} \Biggr).
\end{eqnarray*}
\end{lemma}

We now let $W_n$ be the standardized $ \sum_{k=1}^n\hat{\xi}_k$
given by
%
%e4.38 #&#
%
\begin{equation}\label{stan-sum}
W_n = \frac{1}{\hat{\sigma}_n} \Biggl(\sum_{k=1}^n
\hat{\bolds{\xi }}_k -\hat {\mu}_n \Biggr),
\end{equation}
where $\hat{\mu}_n$ and $\hat{\sigma}_n$ are defined in (\ref{mn}).
As a direct consequence of Lemma \ref{rci} by letting $\omega_k=(\hat
{\xi}_k - \e\hat{\xi}_k)/\hat{\sigma}_n$,
\[
\Delta_1= -\hat{\mu}_n/\hat{\sigma}_n +
y_n \hat{S}_{n,1}\hat {S}_{n,2}/(n \hat{
\sigma}_n) - 2n^{\beta}y_n^3
\varepsilon_{n1}^2/\hat{\sigma}_n,\qquad
\Delta_2= -\hat{\mu}_n/\hat{\sigma}_n
\]
and
\[
\hat{S}_{n,1}^{(k)} = \hat{S}_{n,1}-
\hat{x}_{k1},\qquad
\hat {S}_{n,2}^{(k)} =
\hat{S}_{n,2} - \hat{x}_{k2}^{\tau},\qquad
1\leq k\leq
n,
\]
we have
\begin{eqnarray*}
&&
P \Biggl\{ (y_n/n)\hat{S}_{n,1}
\hat{S}_{n,2} - 2n^{\beta
}y_n^3
\varepsilon_{n1}^2 \leq\sum_{k=1}^n
\hat{\xi}_k <0 \Biggr\}
\\
&&\qquad\leq 21 \Biggl( \hat{\sigma}_n^{-3} \sum
_{k=1}^n\e|\hat{\xi }_k|^3
+ y_n(n\hat{\sigma}_n)^{-1} \e|
\hat{S}_{n,1}\hat{S}_{n,2}|\\
&&\hspace*{16.7pt}\qquad\quad{}  + (\log p)^2n^{-3/2}
+ y_n n^{-1}\hat{\sigma}_n^{-2}
\sum_{k=1}^n \e\bigl|\hat{\xi}_k
\hat{x}_{k1}\hat{S}_{n,2}^{(k)} + \hat{
\xi}_k \hat{x}_{k2}^{\tau}\hat{S}_{n,1}^{(k)}\bigr|\\
&&\qquad\quad\hspace*{151.5pt}{} + y_n n^{-1}\hat{\sigma}_n^{-2}
\sum_{k=1}^n\e\bigl|\hat{\xi}_k
\hat {x}_{k1}\hat {x}_{k2}^{\tau}\bigr| \Biggr)
\\
&&\qquad\leq C \Biggl( n^{-1/2}+ n^{-3/2} \bigl( \e
\hat{S}_{n,1}^2 \bigr)^{1/2} \cdot \bigl( \e\hat
{S}_{n,2}^2 \bigr)^{1/2}
\\
&&\hspace*{13pt}\qquad\quad{} + n^{-2} \sum_{k=1}^n \bigl
\{ \e\hat{S}_{n,1}^{(k) 2} \bigr\}^{1/2}+
n^{-2} \sum_{k=1}^n \bigl\{ \e
\hat{S}_{n,2}^{(k) 2} \bigr\}^{1/2} \Biggr)
\\
&&\qquad\leq Cn^{-1/2}.
\end{eqnarray*}
This, together with expressions (\ref{kn2}) and (\ref{kn21}), verify our
claim (\ref{kn-asym}) and thus
complete the proof of case 2.
\end{pf*}

%s5 #&#
\section{\texorpdfstring{Proof of Theorem \protect\ref{thm3}}{Proof of Theorem 2.3}}\label{sec5}

The main idea of the proof is similar to that of Theorem \ref{thm2}.
%% However, since the normality assumption is removed, necessary
%modifications are required
%%based on results from Theorem \ref{thm2}.
We start with the following three technical lemmas, and their proofs are
%%which are improved analogue of Lemmas 6.9, 6.10 and 6.11 in
postponed to the end of this section.

Let $\{ (z_{k1}, z_{k2}, z_{k3}, z_{k4})^T; k\geq1 \}$ be a sequence
of i.i.d. random vectors with mean zero and common covariance matrix
$\Sigma_4$, which will be specified under different settings. Set
\[
D_{n,i}^2 =\sum_{k=1}^n
z_{ki}^2,\qquad i\in\{1, 2, 3, 4\}.
\]
Suppose $p=p_n\rightarrow\infty$, $\log p=o(n^{\beta})$ as $n
\rightarrow\infty$. For $y\in\mathbb{R}$, let
%
%e5.1 #&#
%
\begin{equation}\label{yn}
y_n = \cases{ \sqrt{(y+4\log p -\log_2 p) /n}, &\quad $0<
\alpha\leq1$, \vspace*{2pt}
\cr
\sqrt{(y+4\log p+c_{n,p}-
\log_2 p) /n}, &\quad $1< \alpha\leq4/3$,}
\end{equation}
for large $n$, where $c_{n,p}=(8\kappa^2/3)n^{-1/2}(\log p)^{3/2}$.

%le5.1 #&#
\begin{lemma} \label{69}
Assume
\[
\Sigma_4= \pmatrix{ 1 & 0 & r & 0
\cr
0 & 1 & 0 & 0
\cr
r & 0 & 1 &
0
\cr
0 & 0 & 0 & 1},\qquad |r| \leq1.
\]
Then, for any $0< \varepsilon<1$,
\[
\sup_{|r|\leq1} P \biggl( \frac{|\sum_{k=1}^n
z_{k1}z_{k2}|}{D_{n,1}D_{n,2}} > y_n,
\frac{|\sum_{k=1}^n z_{k3}z_{k4}|}{D_{n,3}D_{n,4}} > y_n \biggr) = O\bigl(p^{-4(1-\varepsilon)}\bigr).
\]
%
%%as $n\rightarrow\infty$.
\end{lemma}

%le5.2 #&#
\begin{lemma} \label{610}
Assume
\[
\Sigma_4= \pmatrix{ 1 & 0 & r_1 & 0
\cr
0 & 1 &
r_2 & 0
\cr
r_1 & r_2 & 1 & 0
\cr
0 & 0 & 0
& 1},\qquad |r_1| \leq1,\qquad |r_2|\leq1.
\]
Then, for any $0<\varepsilon<1$,
\[
\sup_{|r_1|, |r_2|\leq1} P \biggl( \frac{|\sum_{k=1}^n
z_{k1}z_{k2}|}{D_{n,1}D_{n,2}} > y_n,
\frac{|\sum_{k=1}^n z_{k3}z_{k4}|}{D_{n,3}D_{n,4}} > y_n \biggr) = O\bigl(p^{-4(1-\varepsilon)}\bigr).
\]
%
%%as $n\rightarrow\infty$.
\end{lemma}

%le5.3 #&#
\begin{lemma} \label{611}
Assume
\[
\Sigma_4= \pmatrix{ 1 & 0 & r_1 & 0
\cr
0 & 1 & 0 &
r_2
\cr
r_1 & 0 & 1 & 0
\cr
0 & r_2 & 0 &
1},\qquad |r_1| \leq1,\qquad |r_2|\leq1.
\]
Then, for any $\delta\in(0, 1)$,
\[
\sup_{|r_1|, |r_2| \leq1-\delta} P \biggl( \frac{|\sum_{k=1}^n
z_{k1}z_{k2}|}{D_{n,1}D_{n,2}} > y_n,
\frac{|\sum_{k=1}^n z_{k3}z_{k4}|}{D_{n,3}D_{n,4}} > y_n \biggr) = O\bigl(p^{-2(1+\varepsilon_{\delta})}\bigr),
\]
where
\[
\varepsilon_{\delta}= \bigl( 2 \delta- \delta^2\bigr)/\bigl(4-2
\delta+\delta^2\bigr).
\]
\end{lemma}

Back to the proof of Theorem \ref{thm3}, w.l.o.g., we assume $\mu=0$
and $\sigma^2=1$.
Following the arguments for Theorem \ref{thm2}, we sketch the proof as follows:

\textit{Step} 1: We have
\[
P \Bigl( \max_{1\leq i< j\leq p, j-i \geq m} |\rho_{ij}| \leq
y_n \Bigr) \rightarrow e^{-e^{-y/2}/{\sqrt{8\pi}}}\qquad \mbox{as } n \rightarrow
\infty.
\]
Set
%
%e5.2 #&#
%
\begin{equation}\label{Lp}
\Lambda_p= \bigl\{ (i,j)\dvtx  1\leq i<j\leq p, j-i \geq m, i, j
\notin\Gamma_{p,\delta} \bigr\}
\end{equation}
and
%
%e5.3 #&#
%
\begin{equation}\label{Ln}
L_n' = \max_{(i,j) \in\Lambda_p} |
\rho_{ij}|.
\end{equation}
Clearly,
%
%e5.4 #&#
%
\begin{eqnarray}\label{connect}
P\bigl(L_n' > y_n\bigr) &\leq& P \Bigl(
\max_{1\leq i< j\leq p, j-i \geq m} |\rho_{ij}| > y_n \Bigr)
\nonumber\\[-8pt]\\[-8pt]
&\leq& P\bigl(L_n' > y_n\bigr) + \sum
P \bigl( |\rho_{ij}| >y_n \bigr),\nonumber
\end{eqnarray}
where the last summation is carried out over all pairs $(i,j)$ such
that $1\leq i<j\leq p, j-i\geq m$
and either $i$ or $j$ is in $\Gamma_{p,\delta}$. The total number of
such pairs is no more
than $2p\llvert \Gamma_{p,\delta}\rrvert =o(p^2)$.

Under $H_0$, $\mathbf{x}_1$ and $\mathbf{x}_{m+1}$ are independent
and identically distributed.
Then, by (\ref{PA12}) and (\ref{newp12}), we have for all $0<\alpha
\leq4/3$,
%
%e5.5 #&#
%
\begin{equation}\label{1tau}
P \bigl( |\rho_{1,m+1}| >y_n \bigr) \thicksim\frac{e^{-y/2}}{\sqrt{2\pi}}p^{-2},
\end{equation}
which, in turn, implies that the last summation in (\ref{connect}) is $o(1)$.

\textit{Step} 2: In view of (\ref{connect}) and (\ref{1tau}), it suffices
to prove
%
%e5.6 #&#
%
\begin{equation}\label{step2}
P\bigl(L_n' \leq y_n\bigr) \rightarrow
e^{-e^{-y/2}/{\sqrt{8\pi}}}.
\end{equation}
We follow the lines of proof of Proposition 6.4 in \citet{CJ11a} with
the help of Lemma \ref{pa} and Lemmas \ref{69}--\ref{611}. For
$(i,j)\in\Lambda_p$, set
\[
B_{i,j} = \bigl\{ (k, l) \in\Lambda_p \setminus\bigl\{(i,
j)\bigr\}; \min\bigl\{ |k-i|, |l-j|, |k-j|, |l-i|\bigr\}< m \bigr\}
\]
and $A_{ij}=\{|\rho_{ij}| > y_n\}$ with $y_n$ given in (\ref{yn}). Note
that $|B_{i,j}|\leq4 \times(2m \times p)=8m p$ and $(\mathbf{x}_i,
\mathbf{x}_j)$ are independent of $\{(\mathbf{x}_k, \mathbf{x}_l);
(k, l) \in\Lambda_p \setminus B_{i,j} \}$. By Lemma \ref{pa},
%
%e5.7 #&#
%
\begin{equation}\label{pa3}
\bigl| P\bigl(L_n' \leq y_n \bigr) -
e^{-\lambda_n} \bigr| \leq b_{n,1} +b_{n,2},
\end{equation}
where
%
%e5.8 #&#
%
\begin{eqnarray}\label{est1}
\lambda_n &=& \llvert \Lambda_p \rrvert P(
A_{1,m+1}),\nonumber\\[-8pt]\\[-8pt]
b_{n,1} &=&
\mathop{\sum_{(i, j)\in\Lambda_p}}_{(k, l)\in B_{i, j}}P(
A_{1,m+1})^2 \leq4m p^3 P( A_{1,m+1})^2\nonumber
\end{eqnarray}
and
%
%e5.9 #&#
%
\begin{equation}\label{est2}
b_{n,2} = \sum_{(i, j)\in\Lambda_p} \sum
_{(k, l)\in B_{i, j}} P ( A_{ij}A_{kl} ).
\end{equation}
Clearly, $\llvert  \{(i,j)\dvtx  j\geq i+m \} \rrvert =(p-m)(p-m+1)/2$ and by
definition (\ref{Lp}),
\[
\bigl| \llvert \Lambda_p \rrvert - \bigl\llvert \bigl\{(i,j)\dvtx  j\geq i+m
\bigr\}\bigr\rrvert \bigr| \leq2p \llvert \Gamma_{p,\delta} \rrvert =o
\bigl(p^2\bigr).
\]
This implies $\llvert  \Lambda_p \rrvert  \sim p^2/2$ by assumption on
$m$, which, together with (\ref{1tau}) gives
%
%e5.10 #&#
%
\begin{equation}\label{est1d}
\lambda_n \sim e^{-y/2}/\sqrt{8\pi} \quad\mbox{and}\quad
b_{n,1} = o(1) \qquad\mbox{as } n\rightarrow\infty.
\end{equation}

It remains to estimate $b_{n,2}$. Fix $(i, j)\in\Lambda_p$ and $(k,
l) \in B_{i,j}$ with $i<j$ and $k<l$.
Without loss of generality, assume $i\leq k$ (the case $k<i$ can be
identically proved),
then by definition of $B_{i,j}$
%
%e5.11 #&#
%
\begin{equation}\label{restric}
\min\bigl\{ k-i, |k-j|, |l-j| \bigr\} < m.
\end{equation}
Consider three different cases for the locations of $(i, j)$ and $(k,
l)$ from the above restrictions:
\begin{longlist}[(3)]
\item[(1)] $i< j \leq k< l$, $k-j< m$;

\item[(2)] $i \leq k<l \leq j$, $\min\{ k-i, j-l \}< m$;

\item[(3)] $i \leq k \leq j \leq l$, $\min\{ k-i, j-k, l-j \}< m$.
\end{longlist}

Let $\Omega_{\nu}$ be the subset of index $(i, j,k,l)$ with
restriction $(\nu)$ for $\nu=1, 2, 3$ and
formulate the estimation of $P(A_{ij} A_{kl})$ into three different
cases accordingly.

\textit{Case} (1). It is easy to see that $\llvert  \Omega_1 \rrvert
\leq m p^3= o(p^{3+\varepsilon_\delta})$.
For fixed $(i, j, k, l) \in\Omega_1$, the covariance matrix of
$(x_{1j}, x_{1i}, x_{1k}, x_{1l})$ is equal to
\[
\pmatrix{ 1 & 0 & r & 0
\cr
0 & 1 & 0 & 0
\cr
r & 0 & 1 & 0
\cr
0 & 0 & 0 & 1}
\]
for some $|r|\leq1$. Now we apply Lemma \ref{69} to bound
$P(A_{ij}A_{kl})$. Put
\[
\hat{\rho}_{st} = \frac{\sum_{k=1}^n
x_{ks}x_{kt}}{V_{n,s}V_{n,t}},\qquad
1\leq s<t\leq p,
\]
and analogously to (\ref{event1}), let
%
%e5.12 #&#
%
\begin{equation}\label{event2}
\mathcal{E}_{n\cdot4} = \Bigl\{ \max_{s \in\{i, j, k,l\}} |
\Delta_{n,s}| \leq\varepsilon_{n2} \Bigr\},
\end{equation}
where $\varepsilon_{n2}$ are chosen of the same type as in (\ref{en12})
such that
$ P(\mathcal{E}_{n\cdot4}^c ) =o(p^{-4})$. On $\mathcal{E}_{n \cdot
4}$, we have
\[
|\rho_{st}| \leq \bigl( |\hat{\rho}_{st}|+
\varepsilon_{n2}^2 \bigr)/\bigl(1-\varepsilon_{n2}^2
\bigr) \qquad\mbox{with } \varepsilon_{n2}^2 \asymp(\log p)/n,
\]
which, together with Lemma \ref{69} and the fact that $y_n \sim
2n^{-1/2}(\log p)^{1/2}$, implies that,
for any $0<\varepsilon<(1-\varepsilon_\delta)/4$ and all sufficiently
large $n$,
%
%e5.13 #&#
%
\begin{eqnarray}\label{joint1}
&&P(A_{ij}A_{kl})
\nonumber
\\
&&\qquad\leq P \bigl( |\hat{\rho}_{ij}|> \bigl(1+o(1)\bigr) y_n,
|\hat{\rho}_{kl}|> \bigl(1+o(1)\bigr)y_n \bigr) +o
\bigl(p^{-4}\bigr)
\\
&&\qquad\leq Cp^{-4(1-\varepsilon)}
\nonumber
\end{eqnarray}
and hence
%
%e5.14 #&#
%
\begin{equation}\label{c1}
\sum_{\Omega_1 }P(A_{ij}A_{kl})
=o(1).
\end{equation}
We remark that the $o(1)$'s appeared in (\ref{joint1}) are of order
$n^{-1/2}(\log p)^{1/2}$.\vspace*{9pt}

\textit{Case} (2). Decompose $\Omega_{2}$ as
\begin{eqnarray*}
\Omega_{2} &=& \bigl\{ (i, j, k, l) \in\Omega_2;
k-i<m, j-l<m \bigr\}
\\
&&{} + \bigl\{ (i, j, k, l) \in\Omega_2; k-i<m, j-l \geq m \bigr\}
\\
&&{} + \bigl\{ (i, j, k, l) \in\Omega_2; k-i \geq m, j-l < m \bigr
\}
\\
&:= & \Omega_{2,a} + \Omega_{2,b}+\Omega_{2,c}.
\end{eqnarray*}
Observe that $\llvert  \Omega_{2,a} \rrvert \leq
m^2p^2=o(p^{2(1+\varepsilon_\delta)})$. For
$(i, j, k, l) \in\Omega_{2,a}$, the covariance matrix of $(x_{1i},
x_{1j}, x_{1k}, x_{1l})$ is equal to
\[
\pmatrix{ 1 & 0 & r_1 & 0
\cr
0 & 1 & 0 & r_2
\cr
r_1 & 0 & 1 & 0
\cr
0 & r_2 & 0 & 1}
\]
for some $|r_1|, |r_2|\leq1-\delta$. Using Lemma \ref{611}, along
the lines of the argument in case (1), we get
\[
P(A_{ij}A_{kl}) \leq Cp^{-2(1+\varepsilon_{\delta})}
\]
and therefore
%
%e5.15 #&#
%
\begin{equation}\label{c21}
\sum_{\Omega_{2,a}} P(A_{ij}A_{kl})
=o(1).
\end{equation}
Clearly, $| \Omega_{2,b} | \leq m p^3$ and $ | \Omega_{2,c} |\leq m
p^3$. For $(i, j, k, l)$ in either
$\Omega_{2,b}$ or $\Omega_{2,c}$, the corresponding covariance matrix
of $(x_{1i}, x_{1j}, x_{1k}, x_{1l})$ is
\[
\mbox{either } \pmatrix{ 1 & 0 & r & 0
\cr
0 & 1 & 0 & 0
\cr
r & 0 & 1 & 0
\cr
0 & 0 & 0 & 1} \quad\mbox{or}\quad \pmatrix{ 1 & 0 & 0 & 0
\cr
0 & 1 & 0 & r
\cr
0 &
0 & 1 & 0
\cr
0 & r & 0 & 1},\qquad |r| \leq1.
\]
By the same argument as that in the proof of (\ref{c1}), we have
%
%e5.16 #&#
%
\begin{equation}\label{c22}
\sum_{\Omega_{2, b} \cup\Omega_{2, c} }P(A_{ij}A_{kl})
=o(1) \qquad\mbox{as } n\rightarrow\infty.
\end{equation}

\textit{Case} (3). We aim to show that
%
%e5.17 #&#
%
\begin{equation}\label{c3}
\sum_{\Omega_3 }P(A_{ij}A_{kl})
=o(1).
\end{equation}
Essentially, this can be done by following similar arguments as in
case (2).
However, for $(i, j, k, l)\in\Omega_3$ which satisfies the restriction
\[
\min\{ k-i, j-k, l-j \}< m,
\]
we need to decompose $\Omega_3$ into seven disjoint subsets and
estimate all
the seven possibilities with the help of Lemmas \ref{69}--\ref{611}
as before. The details are omitted here.

Finally, combining expressions (\ref{c1}), (\ref{c21}), (\ref{c22})
and (\ref{c3}) with (\ref{est2}),
we get $b_{n,2}\rightarrow0$ as $n \rightarrow\infty$. This
completes the
proof of (\ref{step2}).
\end{pf*}

\begin{pf*}{Proof of Lemmas \ref{69}--\ref{611}}
We start with a general consideration for estimating joint
probabilities, and the results in Lemmas \ref{69}--\ref{611} will
follow naturally under various dependence structures. Let
\[
\varepsilon_{n1} = c_1(\log p)^{1/2}/n^{\beta/2}
\]
for some constant $c_1>0$ such that, by (\ref{ldp}),
\[
P \bigl( D_{n,1}^2/n \leq1-\varepsilon_{n1}n^{(\beta-1)/2}
\bigr) =o\bigl(p^{-4}\bigr).
\]
Put $\tilde{y}_n=y_n (1-\varepsilon_{n1}n^{(\beta-1)/2}) \sim2\sqrt {(\log p)/n} $.
Using a similar argument as in the proof of Proposition \ref{prop1}
for estimating $P(A_{12}A_{13})$, we have
%
%e5.18 #&#
%
\begin{eqnarray}\label{dp1}
&&P \biggl( \frac{|\sum_{k=1}^n z_{k1}z_{k2}|}{D_{n,1}D_{n,2}} > y_n, \frac{|\sum_{k=1}^n z_{k3}z_{k4}|}{D_{n,3}D_{n,4}} >
y_n \biggr)
\nonumber
\\
&&\qquad\leq P \biggl(\frac{|\sum_{k=1}^n z_{k1}z_{k2}|}{n} >\tilde{y}_n, \frac{|\sum_{k=1}^n z_{k3}z_{k4}|}{n}
> \tilde{y}_n \biggr) +o\bigl(p^{-4}\bigr)
\nonumber\\[-8pt]\\[-8pt]
&&\qquad\leq P \biggl(\frac{|\sum_{k=1}^n(z_{k1}z_{k2}+z_{k3}z_{k4})|}{n^{1/2}} > 2n^{1/2}\tilde{y}_n
\biggr)
\nonumber\\
&&\qquad\quad{} + P \biggl(\frac{|\sum_{k=1}^n(z_{k1}z_{k2}-
z_{k3}z_{k4})|}{n^{1/2}} > 2n^{1/2}\tilde{y}_n
\biggr) + o\bigl(p^{-4}\bigr).
\nonumber
\end{eqnarray}
Note that $\{z_{k1}z_{k2}+z_{k3}z_{k4}, 1\leq k\leq n\}$ is a sequence
of i.i.d. random variables with mean zero.
\end{pf*}

\begin{pf*}{Proof of Lemmas \ref{69} and \ref{610}}
Under both assumptions on $\Sigma_4$, $z_{14}$ is independent of
$(z_{11}, z_{12}, z_{13})$, so that
\begin{eqnarray*}
\e(z_{11}z_{12}+z_{13}z_{14})^2
&=& \e(z_{11}z_{12})^2+\e (z_{13}z_{14})^2+2
\e[z_{11}z_{12}z_{13}z_{14}]
\\
&=& 2+2\e[z_{11}z_{12}z_{13}] \cdot\e
z_{14} =2.
\end{eqnarray*}
It follows from (\ref{ldp}) that, for any $0<\varepsilon<1$,
\[
P \biggl(\frac{|\sum_{k=1}^n(z_{k1}z_{k2}+z_{k3}z_{k4})|}{n^{1/2}} > 2n^{1/2}\tilde {y}_n \biggr)
\leq2\exp\bigl\{-(1-\varepsilon/2)n\tilde{y}_n^2\bigr\} \leq
2p^{4(1-\varepsilon)}
\]
for all sufficiently large $n$. The second probability in (\ref{dp1})
can be
estimated in exactly the same way, and hence the results of Lemmas \ref
{69} and \ref{610} follow
immediately.
\end{pf*}

\begin{pf*}{Proof of Lemma \ref{611}}
In this case, $(z_{11}, z_{13})$ and $(z_{12}, z_{14})$ are
independent. Then, for all $|r_1|, |r_2|\leq1-\delta$,
\[
\e(z_{11}z_{12}+z_{13}z_{14})^2
=2+2\e[z_{11}z_{13}] \cdot\e [z_{12}z_{14}]
\leq2+2(1-\delta)^2.
\]
Set $\varepsilon_{\delta}=(2\delta-\delta^2)/(4-2\delta+\delta^2)$.
Applying (\ref{ldp}) again, we have
\begin{eqnarray*}
&&
P \biggl(\frac{|\sum_{k=1}^n(z_{k1}z_{k2}+z_{k3}z_{k4})|}{n^{1/2}} > 2n^{1/2}\tilde{y}_n
\biggr)
\\
&&\qquad \leq 2\exp \biggl\{-\frac{(1-\varepsilon_{\delta}/2)n\tilde
{y}_n^2}{1+(1-\delta)^2} \biggr\} \leq2p^{-{4(1-\varepsilon
_{\delta
})}/({1+(1-\delta)^2})}\\
&&\qquad=2p^{-2(1+\varepsilon_{\delta})}
\end{eqnarray*}
for all sufficiently large $n$. This completes the proof.
\end{pf*}


%suskaldyti doi

% imsref loaded by lrinkeviciute, 2013-08-29 15:11:03
% imsref loaded by lrinkeviciute, 2013-08-29 15:20:03

\printaddresses

\end{document}